\documentclass[a4paper,10pt]{article}

\usepackage[T1]{fontenc}
\usepackage[latin1]{inputenc}
\usepackage[english]{babel}
\usepackage{amsmath,amsthm,amssymb}
\usepackage{geometry}
\usepackage[all]{xy}

\title{Galois coverings of weakly shod
  algebras}
\author{Patrick Le Meur}
\date{\today}

\theoremstyle{definition}
\newtheorem{dfn}{Definition}[section]
\theoremstyle{plain}
\newtheorem{Thm}[]{Theorem}
\newtheorem{Cor}[Thm]{Corollary}

\newtheorem{lem}[dfn]{Lemma}
\newtheorem{prop}[dfn]{Proposition}
\theoremstyle{remark}
\newtheorem{ex}[dfn]{Example}
\newtheorem{rem}[dfn]{Remark}

\def\ts#1{\normalfont{\textsf{#1}}}
\def\sq{\null\hfill$\blacksquare$\\}
\makeatletter
\renewcommand\section{\@startsection {section}{1}{0mm}%
{-3.5ex \@plus -1ex \@minus -.2ex}
                                   {0.5ex \@plus.2ex}%
                                   {\normalfont\large\bfseries}}

\makeatother

\begin{document}
\sffamily
\maketitle
\abstract{
We investigate the Galois coverings of weakly shod algebras. For a
weakly shod algebra not quasi-tilted of canonical type, we establish a
correspondence between its Galois coverings and the Galois coverings
of its connecting component. As a consequence we show that a weakly
shod algebra
which is not quasi-tilted of canonical type
is simply connected if and only if its first Hochschild
cohomology group vanishes.}

\section*{Introduction}$\ $

Let $A$ be a finite dimensional $k$-algebra where $k$ is an algebraically closed field.
In order to study the category $\ts{mod}\, A$
of finite dimensional 
(right) $A$-modules we assume that $A$ is basic and
connected. The study of $\ts{mod}\, A$ has risen important classes of
algebras. For example: The \emph{representation-finite} algebras, that is,
with only finitely many isomorphism classes of indecomposable modules; the \emph{hereditary}
algebras, that is, path algebras $kQ$ of finite quivers $Q$ with no
oriented cycles; the \emph{tilted} algebras of type $Q$, that is, endomorphism algebras
$\ts{End}_{kQ}(T)$ of tilting $kQ$-modules (see \cite{happel_ringel});
and the \emph{quasi-tilted} algebras, that is, endomorphism algebras
$\ts{End}_{\mathcal{H}}(T)$ of tilting objects $T$ in a hereditary abelian
category $\mathcal{H}$ (see \cite{happel_reiten_smalo}, a quasi-tilted
algebra which is not tilted is called of canonical type). For the last
three classes, each class is a generalisation of the previous one. More recently, a
new class of algebras has arisen
(see \cite{AC03,reiten_skowronski,S03}): That of laura
algebras. The algebra $A$ is called \emph{laura} if there is an upper bound
in the number of isomorphism classes of indecomposable modules which
can appear in an oriented path of non-zero morphisms between
indecomposable $A$-modules starting 
from an injective  and ending at a projective. It appears
that this class contains the four classes cited above. A laura
algebra which not quasi-tilted is characterised by the existence of a unique
non semi-regular component (that is, containing both a projective and
an injective) in its Auslander-Reiten quiver. It is
called the connecting component as a generalisation of the connecting
components of tilted algebras. Hence a laura algebra
which is not quasi-tilted of canonical type has at least one, and at
most two, connecting components (actually, it
has two if and only if $A$ is concealed). Recall that quasi-tilted
algebras of finite representation type are tilted
(\cite[Cor. 2.3.6]{happel_reiten_smalo}) and those  of infinite
representation type are characterised by the
existence of a sincere separating family of semi-regular standard tubes
(\cite{LS96}). Laura algebras comprise the  \emph{weakly shod} algebras
defined by the existence of an upper bound
for the length of a path of non-zero non-isomorphisms from an 
injective to a projective. Actually  a laura algebra
which is not quasi-tilted is weakly shod if and only
if the connecting
component contains no oriented cycles. Weakly shod algebras
were introduced in  \cite{coelho_lanzilotta} as a generalisation of
shod algebras which were defined in \cite{CL99,RS03} as the class of
algebras for which any indecomposable module has injective dimension
or projective dimension at most $1$. For example, quasi-tilted
algebras are shod and therefore weakly shod.

On the other hand, the covering techniques
(\cite{bongartz_gabriel,riedtmann}) have permitted important progress
in the study of 
representation-finite algebras
(see \cite{bgrs,bretscher_gabriel,gabriel}). These techniques need to
consider algebras as $k$-categories. If  $\mathcal C\to A$
is a Galois covering, then $\ts{mod}\, A$ and $\ts{mod}\, \mathcal C$ are
related by
the so-called push-down functor $F_{\lambda}\colon\ts{mod}\, \mathcal
C\to\ts{mod}\, A$. When $A$ has no
proper Galois covering by a connected and locally bounded $k$-category
(or, equivalently, when the fundamental group of any presentation of
$A$ in the sense of \cite{martinezvilla_delapena} is trivial), we say that $A$ is \emph{simply connected} (see
\cite{assem_skowronski}). Simply connected algebras are of special
interest because of the reduction allowed by the push-down
functors. Also they have been object of many investigations (see 
\cite{assem_skowronski,bautista_larrion_salmeron} for
instance). For example, Bongartz and Gabriel
(\cite{bongartz_gabriel}) have classified the simply connected
representation-finite standard algebras using graded trees. Therefore
a nice characterisation of simply connected algebras would be very
useful. 
In \cite[Pb. 1]{skowronski2}, Skowro\'nski asked the following
question for a tame and triangular algebra $A$: 
\begin{center}
\null  \hfill Is $A$ simply connected if and only if $\ts{HH}^1(A)=0$?\hfill($\mathcal Q$)
\end{center}
Up to now, there have been partial answers to $\mathcal Q$ (regardless the
tame assumption): For algebras derived equivalent to a hereditary algebra in \cite{lemeur6}
(and therefore for
tilted algebras), for tame quasi-tilted algebras in \cite{assem_coelho_trepode}
and for tame weakly shod algebras in
\cite{assem_lanzilotta}. Therefore it is natural to try to answer
$\mathcal Q$ for laura algebras.  This shall be done in a forthcoming text
(\cite{laura}). In the present text we study the case of weakly shod
algebras not quasi-tilted of canonical type, which will serve for the study made in \cite{laura}.
For this purpose we prove the following main result.
\begin{Thm}
\label{thm1}
 Let $A$ be connected, weakly shod and not quasi-tilted of canonical type. Let $\Gamma_A$ be a
 connecting component of $\Gamma(\ts{mod}\, A)$. Let $G$ be a group. Then $A$ admits a
Galois covering with
group $G$ by a connected and locally bounded $k$-category if and only
if $\Gamma_A$ admits a Galois covering with group $G$ of
translation quivers. In particular $A$ admits a Galois
covering with group $\pi_1(\Gamma_A)$ by a connected and locally
bounded $k$-category.
\end{Thm}

By \cite[4.2]{bongartz_gabriel}, the fundamental group $\pi_1(\Gamma_A)$ of a connecting
component $\Gamma_A$ is free and isomorphic to the fundamental group of its
orbit-graph. If $A$ is concealed, then its two connecting components
are the 
unique postprojective and the unique preinjective components, so they
have isomorphic fundamental groups. As a
consequence of our main result we prove that $\mathcal Q$ has a positive
answer for weakly shod algebras.
\begin{Cor}
\label{cor1}
 Let $A$ be connected, weakly
shod and not quasi-tilted of canonical type. Let $\Gamma_A$ be a
connecting component of $\Gamma(\ts{mod}\, A)$. The following conditions are
equivalent:
\begin{enumerate}
\item[(a)] $A$ is simply connected.
\item [(b)] $\ts{HH}^1(A)=0$.
\item [(c)] The orbit-graph $\mathcal{O}(\Gamma_A)$ of $\Gamma_A$
  is a tree.
\item[(d)]  $\Gamma_A$ is simply connected.
\end{enumerate}
\end{Cor}

Our proof of Corollary~\ref{cor1} is independent of the one given in
\cite{assem_lanzilotta} for the tame case. Actually we make no
distinction between the different representation types (finite, tame
or wild). 
The proof of Theorem~\ref{thm1} decomposes in two main steps:
\begin{enumerate}
\item If $F\colon\mathcal C\to A$ is a Galois covering with group $G$, then
  every module $X\in\Gamma_A$ is isomorphic to $F_{\lambda}\widetilde X$
  for some $\widetilde X\in \ts{mod}\, \mathcal C$.  The modules $\widetilde X$, for $X$ in
  $\Gamma_A$, form an
  Auslander-Reiten component of $\mathcal C$. This component is a Galois
  covering with group $G$ of $\Gamma_A$.
\item $A$ admits a Galois covering with group
  $\pi_1(\Gamma_A)$ associated to the universal cover of the
  orbit-graph $\mathcal O(\Gamma_A)$.
\end{enumerate}

As an application of the methods we use, we prove the last main result
of the text.
\begin{Thm}
\label{thm2}
  Let $A'\to A$ be a Galois covering with finite group
  $G$ where $A'$ is a basic and connected finite dimensional
  $k$-algebra . Then:
  \begin{enumerate}
  \item[(a)] $A$ is tilted if and only if $A'$ is tilted.
  \item[(b)] $A$ is quasi-tilted if and only if $A'$ is quasi-tilted.
  \item[(c)] $A$ is weakly shod if and only if $A'$ is weakly shod.
  \end{enumerate}
\end{Thm}

The text is organised as follows. In Section~$1$ we fix some notations
and recall some useful definitions. In Section~$2$ we give some
preliminary results: First, we prove some useful facts on covering
techniques; second, we compare
 the Auslander-Reiten quiver of $A$ and the one of $B$ when
$A=B[M]$. Section~$3$ is the very core of the text and is devoted to
the first step of Theorem~\ref{thm1}. In Section~$4$
 we proceed the second step. In Section~$5$, we prove Theorem~\ref{thm1} and
 Corollary~\ref{cor1}. Finally, we prove Theorem~\ref{thm2} in Section~$6$.

\section{Definitions and notations}
\subparagraph{Notations on $k$-categories}$\ $

We refer the reader to \cite[2.1]{bongartz_gabriel} for notions on $k$-categories
and locally bounded $k$-categories. All locally bounded $k$-categories
are small
and all functors between $k$-categories are $k$-linear (of course, our
module categories will be skeletally small). For a locally bounded
$k$-category $\mathcal C$, its objects set is denoted by $\mathcal C_o$ and the
space of morphisms from an object $x$ to an object $y$ is denoted by
$\mathcal C(x,y)$. If $A$ is a basic finite dimensional $k$-algebra, it is
equivalently a locally bounded $k$-category as follows. Fix  a complete set  
$\{e_1,\ldots,e_n\}$ of pairwise orthogonal
primitive idempotents. Then $A_o=\{e_1,\ldots,e_n\}$
and $A(e_i,e_j)=e_jAe_i$ for every $i,j$. In the sequel, $A$ will
always denote a basic finite dimensional
$k$-algebra.

\subparagraph{Notations on modules}$\ $

Let $\mathcal C$ be a $k$-category. Following \cite[2.2]{bongartz_gabriel},
a (right) $\mathcal C$-module is a $k$-linear
functor
$M\colon\mathcal C^{op}\to \ts{MOD}\, k$ where $\ts{MOD}\, k$ is the category of
$k$-vector spaces. If $\mathcal C'$ is another $k$-category, a $\mathcal C-\mathcal
C'$-bimodule is a $k$-linear functor $\mathcal C\times\mathcal C'^{op}\to
\ts{MOD}\,k$. We write $\ts{MOD}\, \mathcal C$ for the category of $\mathcal
C$-modules and $\ts{mod}\, \mathcal C$ for the full subcategory of finite
dimensional $\mathcal C$-modules, that is, those modules $M$ such that
$\sum\limits_{x\in \mathcal C_o} \ts{dim}_k\,
M(x)<\infty$.  The standard duality $\ts{Hom}_k(-,k)$ is denoted by $D$. We
write $\ts{ind}\,\mathcal C$ for the full subcategory of $\ts{mod}\, \mathcal
C$ containing exactly one representative of each isomorphism class of
indecomposable modules. 
A set $\mathcal X$ of modules is called \textit{faithful} if
$\bigcap\limits_{X\in \mathcal X}\ts{Ann}(X)=0$
where $\ts{Ann}(X)$ is the annihilator of $X$, that is, the $\mathcal C-\mathcal
C$-subbimodule of $\mathcal C$ such that $\ts{Ann}(X)(x,y)=\{u\in\mathcal C(x,y)\
|\ mu=0\ \text{for every $m\in X(y)$}\}$.
If $S$ is a set of finite dimensional $\mathcal C$-modules, then
$\ts{add}(S)$ denotes the smallest full subcategory of $\ts{mod}\, \mathcal C$ containing
$S$ and stable under direct sums and direct summands.

Assume that $\mathcal C$ is locally bounded. We write $\Gamma(\ts{mod}\, \mathcal C)$
for the Auslander-Reiten quiver and $\tau_{\mathcal C}=D\ts{Tr}$ for the
Auslander-Reiten translation. Let $\Gamma$ be a component of $\Gamma(\ts{mod}\, \mathcal
C)$. Then $\Gamma$  is called 
\textit{generalised standard} if 
 $\ts{rad}^{\infty}(X,Y)=0$ for every $X,Y\in\Gamma$ (see
 \cite{skowronski3}). 
Here $\ts{rad}$ denotes the radical of $\ts{mod}\,\mathcal C$, that is, the
ideal generated by the non-isomorphisms between indecomposable
modules, $\ts{rad}^n$ denotes the $n$-th power of the radical and
$\ts{rad}^{\infty}=\bigcap\limits_{n\geqslant 1}\ts{rad}^n$.
The component $\Gamma$ is called \textit{non semi-regular} if it
contains both an injective module and a projective module.
We recall the definition of the \emph{orbit-graph} $\mathcal O(\Gamma)$ of
$\Gamma$ in the case $\Gamma$ has no periodic module (see
\cite[4.2]{bongartz_gabriel} for the general case). First, fix a
\emph{polarisation} in $\Gamma$, that is, for every arrow $\alpha\colon
x\to y$ in $\Gamma$ with $y$ non-projective we fix an arrow
$\sigma(\alpha)\colon \tau_{\mathcal C} y\to x$ in such a way that $\sigma$ induces
a bijection from the set of arrows $x\to y$ to the set of arrows
$\tau_{\mathcal C} y\to x$ (see \cite[1.1]{bongartz_gabriel}). Then $\mathcal
O(\Gamma)$ is the graph whose vertices are the $\tau_{\mathcal C}$-orbits of the
vertices in $\Gamma$ and such that there is an edge $(X)^{\tau_{\mathcal
    C}}-(Y)^{\tau_{\mathcal C}}$ for every $\sigma$-orbit of arrows between
a module in $(X)^{\tau_{\mathcal C}}$ and a module in $(Y)^{\tau_{\mathcal C}}$.

We refer the reader to \cite[Chap. VIII, IX]{ass} for a background on
tilting theory.

\subparagraph{Weakly shod algebras (\cite{coelho_lanzilotta})}$\ $

 Let $\mathcal C$ be a locally
bounded $k$-category and $X,Y\in \ts{ind}\,\mathcal C$. A path
$X\rightsquigarrow Y$ in $\ts{ind}\,\mathcal C$ (or in $\Gamma(\ts{mod}\, \mathcal C)$) is a
sequence of non-zero morphisms (or of irreducible morphisms, respectively)
between indecomposable $\mathcal C$-modules
$X=X_0\xrightarrow{f_1}X_1\to\ldots \to X_{n-1}\xrightarrow{f_n}
  X_n=Y$ (with $n\geqslant 0$).
We then say that \textit{$X$ is a predecessor of
  $Y$} and that \textit{$Y$ is a successor of $X$} in $\ts{ind}\,\mathcal C$
(or in $\Gamma(\ts{mod}\, \mathcal C)$, respectively). Hence $X$ is a
successor and a predecessor of itself.

 The
algebra $A$ is called \textit{weakly shod} if
and only if the length of paths in $\ts{ind}\,A$ from an injective to a
projective is bounded. We write $\mathcal P_A^f$ for the set of indecomposable projectives
 which are successors of indecomposable injectives. When $A$ is weakly
 shod this set is
 partially ordered (\cite[4.3]{assem_lanzilotta}) by the relation:
   $P\leqslant Q$ if and only if $P$ is a predecessor of $Q$ in $\ts{ind}\,A$.
We need the following properties when $A$ is weakly shod and connected:
\begin{enumerate}
\item[(a)] If $\mathcal P_A^f=\emptyset$, then $A$ is quasi-tilted
  (\cite[Thm. II.1.14]{happel_reiten_smalo}).
\item[(b)] If $\mathcal P_A^f\neq\emptyset$, then $\Gamma(\ts{mod}\, A)$ has a unique non
  semi-regular component
  (\cite[1.6, 5.4]{coelho_lanzilotta}). This component is generalised standard,
   faithful, has no oriented cycle and contains all the modules lying
   on a path in $\ts{ind}\,A$ form an injective to a projective. In
   particular, every module
  lying on it is a brick (\cite[IV.1.4]{ass}). This
  component is called \textit{the connecting component of
    $\Gamma(\ts{mod}\, A)$ (or of $A$).}
\end{enumerate}
Assume that $A$ is connected, weakly shod and that  $\mathcal
P_A^f\neq \emptyset$. Let $P_m\in\mathcal P_A^f$ be maximal and $e$
the idempotent such that $P_m=eA$. Then $A=B[M]$ where $M=\ts{rad}(P_m)$
and $B=(1-e)A(1-e)$. Moreover:
\begin{enumerate}
\item[(a)] Any component $B'$ of $B$ is weakly
  shod. It is moreover tilted if $\mathcal P_{B'}^f=\emptyset$ (\cite[4.8]{coelho_lanzilotta}).
\item[(b)] Let $M'\in \ts{ind}\,B$ be a summand of $M$ and
  $B'$ the component of $B$ such that $M'\in
  \ts{ind}\,B'$. Then $B'$ is weakly shod
  and not quasi-tilted of canonical type and $M'$ lies on a
  connecting component of $\Gamma(\ts{mod}\, B')$ (\cite[5.3]{assem_lanzilotta}).
\end{enumerate}
Recall (\cite[Thm. 3.1]{reiten_skowronski}) that if a connected algebra $A$ admits a non semi-regular
component which is faithful, generalised standard and has no oriented
cycle, then $A$ is weakly shod.

\subparagraph{Galois coverings of translation quivers
    (\cite{bongartz_gabriel,riedtmann})}$\ $

Let $\Gamma$ and $\Gamma'$ be translation quivers and assume that
$\Gamma$ is connected. A \textit{covering of translation quivers}
$p\colon \Gamma'\to \Gamma$ is a morphism of quivers such that: (a)
$p$ is a covering of unoriented graphs; (b) $p(x)$ is projective
(or injective, respectively) if and only if so is $x$; (c) $p$
commutes with the translations. It is called a \textit{Galois covering
with group $G$} if, moreover, the group $G$ acts on $\Gamma'$ in such a way that:
(d) $G$ acts freely on vertices; (e) $p\,g=p$ for every $g\in G$;
(f) the translation quiver morphism $\Gamma'/G\to\Gamma$ induced by
$p$ is an isomorphism; (g) $\Gamma'$ is connected. Given a connected
translation quiver $\Gamma$, there exists a group
$\pi_1(\Gamma)$ (called the \textit{fundamental group of $\Gamma$}) and a Galois covering $\widetilde{\Gamma}\to
\Gamma$ with group $\pi_1(\Gamma)$ called the \textit{universal cover of
  $\Gamma$}, which factors through any covering $\Gamma'\to
\Gamma$. If $p\colon \Gamma'\to \Gamma$ is a covering (or a
Galois covering with group $G$), then
it naturally induces a covering (or a Galois covering with group $G$, respectively)
$\mathcal O(\Gamma')\to \mathcal O(\Gamma)$ between the associated orbit-graphs. It is
proved in \cite[4.2]{bongartz_gabriel} that if $\Gamma$ has only finitely
many $\tau$-orbits and if $p\colon\Gamma'\to \Gamma$ is the universal
cover of translation quiver, then $\mathcal O(\Gamma')\to \mathcal O(\Gamma)$ is the
universal cover of graphs, that is, $\pi_1(\Gamma)$ is isomorphic to
$\pi_1(\mathcal O(\Gamma))$ (and therefore is free).

\subparagraph{Group actions on module categories (\cite{gabriel})}$\ $

Let $G$ be a group. A \textit{$G$-category} is a $k$-category $\mathcal C$
together with a group morphism $G\to \text{\normalfont\textsf{Aut}}(\mathcal C)$. This defines an
action of $G$ on $\ts{MOD}\, \mathcal C$: If $M\in \ts{MOD}\, \mathcal C$ and $g\in G$, then
$^gM=M\circ g^{-1}$. We write $G_M:=\{g\in G\ |\ \,^gM\simeq M\}$ for
the \textit{stabiliser} of $M$. We say that \emph{$G$ acts freely on
  $\mathcal C$} if the induced action on $\mathcal C_o$ is free.
Assume that $\mathcal C$ is locally bounded. Then this
$G$-action preserves Auslander-Reiten sequences and commutes with
$\tau_{\mathcal C}$. Also it induces an action on $\Gamma(\ts{mod}\, \mathcal
C)$ and on $\mathcal{O}(\Gamma)$ for any $G$-stable component $\Gamma$ of
$\Gamma(\ts{mod}\, \mathcal C)$.

\subparagraph{Galois coverings of categories (\cite{gabriel})}$\ $

Let $G$ be a group and $F\colon\mathcal E\to\mathcal B$ a functor
between $k$-categories. We set
  $\ts{Aut}(F)=\{g\in \ts{Aut}(\mathcal E)\ |\ F\circ g=F\}$. We say that
  $F$ is a 
  \textit{Galois covering with group $G$} if there is a group
  morphism $G\to \text{\normalfont\textsf{Aut}}(F)$ such that $G$ acts
  freely on $\mathcal E$ and
  the induced functor $\overline{F}\colon \mathcal E/G\to \mathcal B$ is
  an isomorphism.  We need the following characterisation for a
  functor $F\colon\mathcal E\to \mathcal B$ to be a Galois covering (\cite[Sect. 3]{gabriel}).
The group morphism $G\to
\text{\normalfont\textsf{Aut}}(F)$ is such that $F$ is Galois with group $G$ if and only if: (a)
the fibres $F^{-1}(x)$ ($x\in \mathcal B_o$) are non-empty and $G$ acts on
these freely and transitively and (b)  $F$ is a \textit{covering
  functor} in the sense of \cite[3.1]{bongartz_gabriel}, that is, for every $x,y\in \mathcal
E_o$ the two maps 
    $\bigoplus\limits_{g\in G}\mathcal E(x,gy)\to\mathcal B(Fx,Fy)$ and 
    $\bigoplus\limits_{g\in G}\mathcal E(gy,x)\to\mathcal B(Fy,Fx)$
    induced by $F$ are isomorphisms.
A Galois covering $F\colon\mathcal E\to\mathcal B$ with $\mathcal E$ and $\mathcal B$ locally
  bounded and connected is called
  \textit{connected}. In such a case, the morphism $G\to \text{\normalfont\textsf{Aut}}(F)$ is an
  isomorphism (\cite[Prop. 6.1.37]{lemeur_thesis}).
A connected and locally bounded $k$-category $\mathcal B$ is called \textit{simply
  connected} if and only if there is no connected Galois covering $\mathcal
E\to\mathcal B$ with non trivial group. This definition is equivalent
(\cite[Cor. 4.5]{lemeur2}) to the
original one of \cite{assem_skowronski}
and it is more convenient for our purposes.

\subparagraph{Covering techniques (\cite{bongartz_gabriel} and
  \cite{gabriel})}$\ $

 Let $F\colon\mathcal E\to \mathcal B$ be a Galois covering
between locally bounded $k$-categories. 
We write $F_{\lambda}\colon \ts{MOD}\, \mathcal E\to \ts{MOD}\, \mathcal B$ and $F_.\colon
\ts{MOD}\, \mathcal B\to \ts{MOD}\, \mathcal E$ for the
\textit{push-down} functor and the 
\textit{pull-up} functor, respectively. Recall
(\cite[3.2]{bongartz_gabriel}) that $F_.=X\circ F$ for every
$X\in\ts{MOD}\,\mathcal B$ and that for $M\in\ts{MOD}\,\mathcal E$, the $\mathcal B$-module
$F_{\lambda}M$ is such that $F_{\lambda}M(x)=\bigoplus\limits_{Fx'=x}M(x')$ for
every $x\in\mathcal B_o$.
We list some needed properties on theses functors. Both $F_{\lambda}$ and $F_.$
are exact; $(F_{\lambda},F_.)$ is adjoint; $F_{\lambda}M$ is projective (or
injective) if and only if $M$ is projective (or injective,
respectively); $F_{\lambda}(\ts{mod}\,\mathcal E)\subseteq \ts{mod}\,\mathcal B$; the
functor $F_{\lambda}$ is $G$-invariant, that is, $F_{\lambda}\circ g=F_{\lambda}$ for
every $g\in G$; for every $M\in\ts{mod}\,\mathcal E$ we have
$F_.F_{\lambda}M\simeq\bigoplus\limits_{g\in G}\,^gM$
(\cite[3.2]{gabriel}); and $F_{\lambda}$ commutes with the duality,
that is, 
$D\circ F_{\lambda}\simeq F^{op}_{\lambda}\circ D$ on $\ts{mod}\,\mathcal E$. Finally, it satisfies a property
  which will be refered to as the \textit{covering property of $F_{\lambda}$}: For
  $M,N\in \ts{mod}\,\mathcal E$, the two maps $\bigoplus\limits_{g\in
    G}\ts{Hom}_{\mathcal E}(\,^gM,N)\to \ts{Hom}_{\mathcal B}(F_{\lambda}M,F_{\lambda}N)$ and
      $\bigoplus\limits_{g\in G}\ts{Hom}_{\mathcal E}(M,\,^gN)\to
      \ts{Hom}_{\mathcal B}(F_{\lambda}M,F_{\lambda}N)$ induced by 
  $F_{\lambda}$ are $k$-linear isomorphisms.
 A module $X\in
\ts{ind}\,\mathcal B$ is called \textit{of the first kind (with respect to $F$)} if and only
if there exists $\widetilde{X}\in \ts{mod}\,\mathcal E$ (necessarily
indecomposable) such that
$F_{\lambda}\widetilde{X}\simeq X$ in $\ts{mod}\,B$. Note that
if $\widetilde{X}$ exists, then $X=F_{\lambda}\overline{X}$ for some $\overline{X}\in
  \ts{ind}\,\mathcal E$; and,
if $\widetilde{X},\widehat{X}\in \ts{ind}\,\mathcal E$ are such that
  $F_{\lambda}\widetilde{X}\simeq F_{\lambda}\widehat{X}\simeq X$,
  then $\widetilde{X}\simeq\,^g\widehat{X}$ for some $g\in G$ (see \cite[3.5]{gabriel}).

\section{Preliminaries}
\subparagraph{Some results on covering techniques}$\ $

Let $F\colon\mathcal C\to A$ be a Galois covering with group $G$ where
$\mathcal C$ is locally bounded. We prove some useful comparisons between 
of $\Gamma(\ts{mod}\, A)$ and 
$\Gamma(\ts{mod}\, \mathcal C)$. 
First, we give  a necessary condition
on a morphism in $\ts{mod}\,\mathcal C$ to be mapped by $F_{\lambda}$ to a section or
a to retraction.
\begin{lem}
  \label{lem1.1}
Let $X,Y\in \ts{mod}\, \mathcal C$ and $f\in \ts{Hom}_{\mathcal C}(X,Y)$. 
\begin{enumerate}
\item[(a)] $F_{\lambda}(f)$ is a section (or a
retraction) if and only if so is  $f$.
\item[(b)] If $F_{\lambda}(f)$
  is irreducible, then so is $f$.
\item[(c)] Let $u\colon E\to X$ (or $v\colon X\to E$) be a right (or left) minimal
   almost split morphism in $\ts{mod}\, \mathcal C$. Assume that
  $G_X=1$. Then so is $F_{\lambda}(u)$
  (or $F_{\lambda}(v)$, respectively).
\item[(d)] $F_{\lambda}\tau_{\mathcal C}X\simeq \tau_AF_{\lambda}X$.  
\end{enumerate}
\end{lem}
\noindent{\textbf{Proof:}} (a) Obviously, if $f$ is a section
(or a retraction), then so is $F_{\lambda}(f)$. Assume that $F_{\lambda}(f)$
is a section. So $\ts{Id}_{F_{\lambda}X}=r\circ F_{\lambda}(f)$ with
$r\in\ts{Hom}_A(F_{\lambda}X,F_{\lambda}Y)$. Moreover,
$r=\sum\limits_gF_{\lambda}(r_g)$ with $(r_g)_{g\in
  G}\in\bigoplus\limits_{g\in G}\ts{Hom}_{\mathcal C}(Y,\,^gX)$, using the
covering propery of $F_{\lambda}$. Therefore
$\ts{Id}_{F_{\lambda}X}=\sum\limits_gF_{\lambda}(r_g\circ f)$. The covering
property of $F_{\lambda}$ then implies that $\ts{Id}_X=r_1\circ f$, that
is, $f$ is a section. Dually, if $F_{\lambda}(f)$ is a retraction, then so
is $f$.

(b) is a direct consequence of (a).

(c) is due to the proof of
\cite[3.6, (b)]{gabriel}.

(d)  follows from the fact that $F_{\lambda}$ is exact, maps
projective modules to projective modules (in particular, $F_{\lambda}$ maps
a minimal projective resolution in $\ts{mod}\, \mathcal C$ to a minimal projective resolution in
$\ts{mod}\, A$) and commutes with the duality.
\sq

\begin{lem}
  \label{lem1.7}
Let $\Gamma$ be a component of $\Gamma(\ts{mod}\, A)$ made of modules of the
first kind and $\widetilde{\Gamma}$  the full
subquiver of $\Gamma(\ts{mod}\, \mathcal C)$ generated by $\{X\in\Gamma(\ts{mod}\, \mathcal C)\ |\
F_{\lambda}X\in \Gamma\}$. Then:
\begin{enumerate}
\item[(a)] Let $u\colon M\to P$ be a right minimal almost split morphism in
  $\ts{mod}\, \mathcal C$ with $P$ indecomposable projective. Then $F_{\lambda}(u)$
  is right minimal almost split.
\item[(b)] Let $X\in\,\widetilde{\Gamma}$ be non projective. Then $F_{\lambda}$
  transforms any almost split sequence ending at $X$
  into an almost split sequence ending at $F_{\lambda}X$.
\item[(c)] Let $u\in \ts{Hom}_{\mathcal C}(X,Y)$ with $X,Y\in
  \widetilde{\Gamma}$. Then $u$ is irreducible if and only if so is
  $F_{\lambda}(u)$.
\item[(d)]  $\Gamma$ is stable under predecessors and 
  under successors in $\Gamma(\ts{mod}\, \mathcal C)$ and under the
  action of $G$.
\end{enumerate}
\end{lem}
\noindent{\textbf{Proof:}} (a) follows from
\cite[3.2]{bongartz_gabriel}.

(b) Fix an almost split sequence $0\to \tau_{\mathcal
  C}X\xrightarrow{\varphi}E\xrightarrow{\theta}X\to 0$ in $\ts{mod}\, \mathcal C$.
By \ref{lem1.1}, (d), we have  an exact sequence
$0\to
  \tau_A
  F_{\lambda}X\xrightarrow{F_{\lambda}(\varphi)}F_{\lambda}E
  \xrightarrow{F_{\lambda}(\theta)}F_{\lambda}X\to 0$ in $\ts{mod}\, A$.
By \ref{lem1.1}, (a), it does not split. Moreover,
$F_{\lambda}X$ is indecomposable and non-projective. Let $v\colon Z\to F_{\lambda}X$ be right minimal almost
split. We only need  to prove that $v$
factors through $F_{\lambda}\theta$.  Write $v\colon Z\to
F_{\lambda}X$ as
$v=\begin{bmatrix}v_1&\cdots&v_n\end{bmatrix}\colon
  Z_1\oplus\ldots\oplus Z_n\to F_{\lambda}X$
where $Z_1,\ldots,Z_n\in \ts{ind}\,A$. We
prove that each $v_i$ factors through $F_{\lambda}\theta$. We have
$Z_i\in\Gamma$ because $v_i$ is
irreducible.
Therefore $Z_i= F_{\lambda}\widetilde{Z}_i$ for some
$\widetilde{Z}_i\in \ts{mod}\,\,\mathcal C$ indecomposable. So 
$v_i=\sum\limits_gF_{\lambda}(v_{i,g})$ where  $(v_{i,g})_{g\in G}\in
\bigoplus\limits_{g\in G}\ts{Hom}_{\mathcal
  C}(\,^g\widetilde{Z}_i,X)$. Note that
$^g\widetilde{Z}_i\not\simeq X$ for every $g\in G$ because $Z_i\not\simeq F_{\lambda}X$. Thus
$v_{i,g}=\theta\circ w_{i,g}$ for some $w_{i,g}\in \ts{Hom}_{\mathcal
  C}(\,^g\widetilde{Z}_i,E)$ for every $g$.
We may assume that $w_{i,g}=0$ if $v_{i,g}=0$. Then $v_i=F_{\lambda}(\theta)\circ \left(\sum\limits_{g\in
    G}F_{\lambda}(w_{i,g})\right)$ where $\sum\limits_{g\in
    G}F_{\lambda}(w_{i,g})\in \ts{Hom}_A(Z_i,F_{\lambda}X)$ for every $i$. Thus
$v_1,\ldots,v_n$ factor through $F_{\lambda}\theta$. Therefore so
does $v$. This proves (b).

(c) is a direct consequence of (a), (b) and
\ref{lem1.1}.

(d) Clearly, $\widetilde{\Gamma}$ is stable under the action
of $G$. We prove
the stability under
predecessors (the proof for successors is
dual). Let $u\in \ts{Hom}_{\mathcal C}(X,Y)$ be irreducible with $X\in \ts{ind}\,\mathcal C$
and $Y\in\widetilde{\Gamma}$. 
We claim that $F_{\lambda}X\in \ts{add}(\Gamma)$.
If $Y$ is
projective, then $X$ is a direct summand of $\ts{rad}(Y)$ and $u\colon X\to
Y$ is the inclusion. So $F_{\lambda}Y$ is indecomposable
projective, $F_{\lambda}X$ is a direct summand of
$F_{\lambda}(\ts{rad}(Y))=\ts{rad}(F_{\lambda}Y)$
(\cite[3.2]{bongartz_gabriel}) and $F_{\lambda}(u)\colon
  F_{\lambda}X\to F_{\lambda}Y$ is injective. Since $F_{\lambda}Y\in
  \Gamma$ we have $\ts{rad}(F_{\lambda}Y)\in \ts{add}(\Gamma)$ and  therefore $F_{\lambda}X\in
  \ts{add}(\Gamma)$.
Assume that $Y$ is not projective. So there is an almost split
sequence in $\ts{mod}\, \mathcal C$:
\begin{equation}
  0\to \tau_{\mathcal C}Y\to E\oplus
  X\xrightarrow{\begin{bmatrix}?\\u\end{bmatrix}}Y\to 0\ .\notag
\end{equation}
By (a), there is an almost split
sequence in $\ts{mod}\, A$:
\begin{equation}
  0\to \tau_A F_{\lambda}Y\to F_{\lambda}E\oplus
  F_{\lambda}X\xrightarrow{\begin{bmatrix}?\\F_{\lambda}u\end{bmatrix}}F_{\lambda}Y\to 0 \ .\notag
\end{equation}
Since $F_{\lambda}Y\in \Gamma$, we have $F_{\lambda}X\in
\ts{add}(\Gamma)$.
This proves the claim. Now we prove that $F_{\lambda}X$ is indecomposable. Since
$F_{\lambda}X\in \ts{add}(\Gamma)$, there exist
$\widetilde{E}_1,\ldots,\widetilde{E}_n\in\widetilde{\Gamma}$ and an
isomorphism $\varphi\colon F_{\lambda}X\xrightarrow{\sim}
\widetilde{E}_1\oplus\ldots\oplus \widetilde{E}_n$. From the covering
property of $F_{\lambda}$, we have $\varphi=\sum\limits_{g\in G}F_{\lambda}(\varphi_g)$
where $(\varphi_g)_{g\in G}\in \bigoplus\limits_{g\in G}\ts{Hom}_{\mathcal
  C}(\,^gX,\widetilde{E}_1\oplus\ldots\oplus
\widetilde{E}_n)$. Since $\varphi$ is an isomorphism, there exists
$g\in G$ such that $F_{\lambda}(\varphi_g)\not\in \ts{rad}(F_{\lambda}X,F_{\lambda}\widetilde
E_1\oplus\ldots\oplus F_{\lambda}\widetilde E_n)$. So there exists
$i$ such that the restriction $F_{\lambda}\,^gX\to F_{\lambda}\widetilde E_i$ of
$F_{\lambda}(\varphi_g)$ is an isomorphism so that $^gX\simeq\widetilde E_i\in\Gamma$.\sq

The following proposition describes the action of $F_{\lambda}$
on almost split sequences in $\ts{mod}\, \mathcal C$ under suitable
conditions. Note that if we assume that $G$ acts freely on
indecomposable $\mathcal C$-modules (that is, $G_X=1$ for any $X\in \ts{ind}\,\mathcal
C$), then the last three points follow at once from \cite[3.6]{gabriel}.
\begin{prop}
  \label{prop1.2}
Keep the hypotheses and notations of \ref{lem1.7}.
\begin{enumerate}
\item[(a)] $\Gamma$ is faithful if and only if $\widetilde{\Gamma}$ is.
\item[(b)] $\Gamma$ is generalised standard if and only if
  $\ts{rad}^{\infty}(X,Y)=0$ for every $X,Y\in \widetilde{\Gamma}$.
\item[(c)] $\widetilde{\Gamma}$ is a
  (disjoint) union of components of $\Gamma(\ts{mod}\, \mathcal C)$. In particular,
  $\widetilde{\Gamma}$ is a translation subquiver of $\Gamma(\ts{mod}\, \mathcal C)$.
\item[(d)] The map $X\mapsto F_{\lambda}X$ extends to a covering of
  translation quivers $\widetilde{\Gamma}\to \Gamma$. If
  $\widetilde{\Gamma}$ is connected and $G_X=1$ for every
  $X\in\widetilde{\Gamma}$, then this is a Galois covering with group $G$.
\item[(e)] $\Gamma$ has an oriented cycle if and only if
  $\widetilde{\Gamma}$ has a non trivial path of the form
  $X\rightsquigarrow\,^gX$.
\end{enumerate}
\end{prop}
\noindent{\textbf{Proof:}} (a) Assume that $\Gamma$ is faithful. Let
$u\in\ts{Ann}(\widetilde{\Gamma})(x,y)$, that is, $u\in\mathcal C(x,y)$ and $mu=0$ for every
$m\in X(y)$, $X\in\widetilde{\Gamma}$. We claim that $F(u)\in\ts{Ann}(\Gamma)(Fx,Fy)$. Let
$X\in\Gamma$ and $m\in X(Fy)$. We may assume that $X=F_{\lambda}\widetilde X$ with $\widetilde
X\in\widetilde{\Gamma}$. So $m=(m_g)_{g\in G}\in\bigoplus\limits_{g\in G}\widetilde X(gy)$
and, therefore, $mF(u)=(m_gg(u))_{g\in G}$. On the other hand,
$g(u)\in\ts{Ann}(\widetilde{\Gamma})(gx,gy)$ because $\widetilde{\Gamma}$ is $G$-stable. So $m_gg(u)=0$
for every $g\in G$ so that $mF(u)=0$. Thus
$F(u)\in\ts{Ann}(\Gamma)(Fx,Fy)=0$ and, therefore, $u=0$. So $\widetilde{\Gamma}$ is
faithful.

Conversely, assume that $\widetilde{\Gamma}$ is faithful and let
$u\in\ts{Ann}(\Gamma)(Fx,Fy)$. So $u=\sum\limits_gF(u_g)$ where
$(u_g)_{g\in G}\in\bigoplus\limits_{g\in G}\mathcal C(gx,y)$. We claim that
$u_g\in\ts{Ann}(\widetilde{\Gamma})(gx,y)$ for every $g\in G$. Indeed, let $X\in\widetilde{\Gamma}$
and $m\in X(y)$. Then $m\in F_{\lambda}X(Fy)$ and $0=mu=(mu_g)_{g\in
  G}\in\bigoplus\limits_{g\in G}F_{\lambda}X(gx)$. So $mu_g=0$ for every
$g$. Thus $u_g\in\ts{Ann}(\widetilde{\Gamma})(gx,y)$ for every $g\in G$ and,
therefore, $u=0$ because $\widetilde{\Gamma}$ is faithful. So $\widetilde{\Gamma}$ is faithful.

(b)  Assume that $\ts{rad}^{\infty}(X,Y)=0$ for every $X,Y\in\widetilde{\Gamma}$. Let
$X,Y\in\widetilde{\Gamma}$. We prove that
$\ts{rad}^{\infty}(F_{\lambda}X,F_{\lambda}Y)=0$. Since
$\ts{Hom}_A(F_{\lambda}X,F_{\lambda}Y)$ is finite dimensional and isomorphic to
$\bigoplus\limits_{g\in G}\ts{Hom}_{\mathcal C}(X,\,^gY)$, there exists
$n\geqslant 1$ such that $\ts{rad}^n(X,\,^gY)=0$ for every $g\in
G$. Let $f\in\ts{rad}^l(F_{\lambda}X,F_{\lambda}Y)$ with $l\geqslant 1$. Let
$\begin{bmatrix}u_1&\ldots &u_t\end{bmatrix}^t\colon X\to
E_1\oplus\ldots\oplus E_t$ be left minimal almost split in $\ts{mod}\,\mathcal
C$. By \ref{lem1.7}, (a) and (b), there exist
$f_i\in\ts{Hom}_A(F_{\lambda}E_i,F_{\lambda}Y)$ for every $i$, such that
$f=\sum\limits_if_i\circ F_{\lambda}(u_i)$. More generally an induction on
$l$ shows that there exist morphisms $\delta_1\colon X\to
X_1,\ldots,\delta_s\colon X\to X_s$ in $\ts{mod}\,\mathcal C$ all equal to
compositions of $l$ irreducible morphisms between indecomposable
modules and there exist $h_i\in\ts{Hom}_A(F_{\lambda}X_i,F_{\lambda}Y)$ for
every $i$, such that $f=\sum\limits_ih_i\circ F_{\lambda}(\delta_i)$. On
the other hand, $h_i=\sum\limits_{g}F_{\lambda}(h_{i,g})$ with
$(h_{i,g})_{g\in G}\in\bigoplus\limits_{g\in G}\ts{Hom}_{\mathcal
  C}(X_i,\,^gY)$ by the covering property of $F_{\lambda}$. Therefore:
\begin{equation}
  f=\sum\limits_gF_{\lambda}\left(\sum\limits_ih_{i,g}\circ\delta_i\right)\notag
\end{equation}
where $\sum\limits_ih_{i,g}\circ\delta_i\in\ts{rad}^l(X,\,^gY)$ for
every $g$. In the particular case where $l=n$, we have $f=0$. Thus
$\ts{rad}^n(F_{\lambda}X,F_{\lambda}Y)=0$. This proves that $\Gamma$ is generalised
standard.

Conversely, assume that $\Gamma$ is generalised standard. Let
$f\in\ts{rad}^l(X,Y)$ with $X,Y\in\widetilde{\Gamma}$ and $l\geqslant 1$. The
arguments used above show that there exist morphisms $\delta_1\colon
X\to X_1,\ldots,\delta_s\colon X\to X_s$ in $\ts{mod}\,\mathcal C$ all equal
to compositions of $l$ irreducible morphisms between indecomposable
modules and there exist morphisms $h_1\colon X_1\to Y,\ldots,h_s\colon
X_s\to Y$ such that $f=\sum\limits_ih_i\circ \delta_i$. By
\ref{lem1.7}, (c), we therefore have
$F_{\lambda}(f)\in\ts{rad}^l(F_{\lambda}X,F_{\lambda}Y)$. Hence
$F_{\lambda}(\ts{rad}^{\infty}(X,Y))\subseteq
\ts{rad}^{\infty}(F_{\lambda}X,F_{\lambda}Y)=0$. Since $F_{\lambda}$ is faithful,
we have $\ts{rad}^{\infty}(X,Y)=0$
for every $X,Y\in\widetilde{\Gamma}$.

(c) This is a direct consequence of \ref{lem1.7}, (d).

(d) By assumption and \ref{lem1.7}, $F_{\lambda}$
preserves indecomposability, irreducibility and almost split sequences
in $\widetilde{\Gamma}$. Consequently, for each $X\in\widetilde{\Gamma}$ there is a
bijection between the set of arrows in $\widetilde{\Gamma}$ which leave (or
arrive at) $X$ and the set of arrows in $\Gamma$ which leave (or
arrive at, respectively) $F_{\lambda}X$. Whence the covering
$\widetilde{\Gamma}\to \Gamma$ extending the map $X\mapsto F_{\lambda}X$.
The rest of the assertion is a consequence of the arguments presented
in the proof of \cite[3.6]{gabriel}.

(e) follows from (d).
\sq

\begin{rem}
  \label{rem1.3}
Assume, in \ref{lem1.7},
that $\widetilde{\Gamma}$ is connected and $G_X=1$ for every $X\in\widetilde{\Gamma}$. By
\ref{prop1.2}, (d), there is a Galois covering with group
$G$ of graphs $p\colon \mathcal O(\widetilde{\Gamma})\to \mathcal O(\Gamma)$ such that
$p\left((X)^{\tau_{\mathcal C}}\right)=\left(F_{\lambda}X\right)^{\tau_A}$ for
every vertex $X\in\widetilde{\Gamma}$. The
$G$-action on $\mathcal O(\widetilde{\Gamma})$ is given by $^g\left((X)^{\tau_{\mathcal
      C}}\right)=\left(\,^gX\right)^{\tau_{\mathcal C}}$ for every $g\in G$,
$X\in\widetilde{\Gamma}$. In particular, if $g\colon\mathcal O(\widetilde{\Gamma})\to\mathcal O(\widetilde{\Gamma})$ is an
automorphism of graphs such that $p\circ g=p$, then there exists
$g'\in G$ such that $g$ is induced by $g'$.
\end{rem}

\begin{rem}
  In view of the proof of \ref{prop1.2}, (a), if $X\in\ts{mod}\,\mathcal C$ is faithful,
  then so is $F_{\lambda}X$. However, one can easily find examples where
  $F_{\lambda}X$ is faithful and $X$ is not.
\end{rem}

\subparagraph{Comparisons between the Auslander-Rieten quivers of $A$ and
  $B$ when $A=B[M]$}$\ $

In this paragraph we assume that $A$ is connected and weakly shod and
$\mathcal P_A^f\neq \emptyset$. Let $P_m\in\mathcal P_A^f$ be
maximal and $A=B[M]$  the associated one-point extension. We
give a useful relationship between the connecting component $\Gamma_A$
of $\Gamma(\ts{mod}\, A)$ and the connecting components associated to the connected
components of $B$. It follows from 
the work made in \cite{coelho_lanzilotta} (see also 
\cite[Lem. 4.1]{assem_lanzilotta} who treated the case where the
extension point is
separating). For convenience, we give the
details below. Note that:
\begin{enumerate}
\item[(a)] Any strict predecessor of $P_m$ in $\ts{ind}\,A$ is a $B$-module.
\item[(b)] If $P\in \ts{ind}\,B$ is projective, then any predecessor of $P$ in
  $\ts{ind}\,A$ is a $B$-module.
\end{enumerate}
We begin with the following lemma.
\begin{lem}
\label{lem1.4}
Let $\mathcal X$ be the full subcategory of $\ts{ind}\,A$ generated by
\begin{equation}
  \{X\in \ts{ind}\,A\ |\ X\not\simeq P_m\ \text{and $X$ is a
    predecessor in $\ts{ind}\,A$ of an indecomposable projective $A$-module}\}\ .\notag
\end{equation}
Then:
\begin{enumerate}
\item[(a)] $\mathcal X$ is made of $B$-modules.
\item[(b)] $\mathcal X$ is stable under predecessors in $\ts{ind}\, A$ and
  contains no successor of $P_m$ in $\ts{ind}\, A$.
\item[(c)] $\tau_A$ and $\tau_B$ coincide on $\mathcal X$.
\item[(d)] The full subquivers of $\Gamma(\ts{mod}\, A)$ and $\Gamma(\ts{mod}\, B)$
  generated by $\mathcal X$ coincide.
\end{enumerate}
\end{lem}
\noindent{\textbf{Proof:}} (a) and (b) follow from the definition of
$\mathcal X$.
For $L\in \ts{mod}\, A$ let $\overline{L}$ be the
$B$-module obtained by restriction of scalars, that is,
$\overline{L}=L.(1-e)$ if $e\in A$ is the idempotent such that $P_m=eA$.
Assume that $0\to \tau_AX\xrightarrow{u}
  E\xrightarrow{v} X\to 0$ is an almost split sequence in
  $\ts{mod}\, A$ with $X\in \ts{ind}\,B$. Then it is easily verified that $\overline{\tau_A
    X}=\tau_B X$ and $0\to \tau_BX\xrightarrow{\overline{u}}
  \overline{E}\xrightarrow{\overline{v}} \overline{X}\to 0$ is
  almost split in $\ts{mod}\, B$. Also, if $X$ is not a
  successor of $P_m$, then the two exact sequences coincide. Then (c)
  and (d) follow from these facts.
\sq

The category $\mathcal X$ of the preceding lemma serves to compare connecting
components as follows.
\begin{lem}
\label{lem1.5}
Let $\mathcal X$ be as in the preceding lemma, $M'\in \ts{ind}\,B$
a direct summand of $M$ and  $B'$ the
component of $B$ such that $M'\in \ts{ind}\,B'$. If $\Gamma'$ is the
 component of $\Gamma(\ts{mod}\, B')$ containing $M'$, then:
\begin{enumerate}
\item[(a)] The connecting component $\Gamma_A$ of $\Gamma(\ts{mod}\, A)$ contains every
  module lying on both $\Gamma'$ and $\mathcal X$.
\item[(b)] The full subquivers of $\Gamma_A$ and $\Gamma'$ generated by
  the modules lying on both $\Gamma'$ and $\mathcal X$ coincide.
\item[(c)] Every $\tau_{B'}$-orbit of $\Gamma'$ contains a module lying on
  $\mathcal X$.
\end{enumerate}
\end{lem}
\noindent{\textbf{Proof:}}  
(a) Let $X$ lie on both $\Gamma'$ and $\mathcal X$. By \cite[1.1]{assem_lanzilotta}, $\tau_{B'}^mX$ is a predecessor in $\Gamma(\ts{mod}\, B')$
(and therefore in $\Gamma(\ts{mod}\, A)$, by \ref{lem1.4}) of $M'$ or of a projective $P\in
\ts{ind}\,B'$ for some $m\geqslant 0$. By
\cite[Lem. 5.3]{coelho_lanzilotta}, $P\in\Gamma_A$. So
$\tau_{B'}^mX\in\Gamma_A$. On the
other hand, \ref{lem1.4}, (c), implies that $\tau_{B'}^mX=\tau_A^mX$. So
$X\in\Gamma_A$.

(b) Let $\mathcal X_1$ and $\mathcal X_2$ be the full subquivers of
$\Gamma_A$ and $\Gamma'$, respectively, generated by the modules
lying on both $\mathcal X$ and $\Gamma'$. By (a), $\mathcal X_1$
and $\mathcal X_2$ have the same vertices. Then \ref{lem1.4}, (d), implies
that $\mathcal X_1=\mathcal X_2$.

(c) is obtained using similar arguments as
those used to prove (a).\sq

\begin{rem}
\label{rem1.6}
Using \ref{lem1.5} we get the following description of the orbit-graph
$\mathcal O(\Gamma_A)$. For simplicity, we write $\mathcal
O(\Gamma_A)\backslash\{(P_m)^{\tau_A}\}$ for the full subgraph of $\mathcal
O(\Gamma_A)$ generated by the vertices different from $(P_m)^{\tau_A}$.
\begin{enumerate}
\item[(a)] Let $B'$ be a component of $B$ and $\Gamma_{B'}$ the (unique)
  connecting component of $B'$ containing a direct summand of $M$.
 Then $\mathcal O(\Gamma_{B'})$ is a component of
  $\mathcal O(\Gamma_A)\backslash\{(P_m)^{\tau_A}\}$ and all the components of
  $\mathcal O(\Gamma_A)\backslash\{(P_m)^{\tau_A}\}$ have this form.
\item[(b)]
 If $X$ is an indecomposable direct summand of $M$ with
  multiplicity $d$, then $(X)^{\tau_B}$ lies on exactly one of the
  connected components of $\mathcal O(\Gamma_A)\backslash\{(P_m)\}^{\tau_A}$
  and $\mathcal O(\Gamma_A)$ contains exactly $d$ edges
  $(X)^{\tau_A}-(P_m)^{\tau_A}$. Moreover all the arrows connected to
  $(P_m)^{\tau_A}$ have this form.
\end{enumerate}
\end{rem}

\section{Components of the first kind for weakly shod algebras}$\ $

Let $A$ be weakly shod. We examine when a component of
$\Gamma(\ts{mod}\, A)$ is made of modules of the first kind with respect to any Galois
covering of $A$. We study two cases: When the component is connecting
and when it is semi-regular and not regular.
\subparagraph{Connecting components of the first kind}$\ $

The aim of this paragraph is to prove the following proposition.
\begin{prop}
\label{prop2.1}
  Let $A$ be connected, weakly shod and not quasi-tilted of canonical
  type, $\Gamma_A$  a connecting component of $A$ and
  $F\colon \mathcal C\to A$ a connected Galois covering with group $G$. Then
  $\Gamma_A$ is made of modules of the first kind. Moreover the full
  subquiver $\Gamma_{\mathcal C}$ of $\Gamma(\ts{mod}\,\mathcal C)$ generated by the
  modules $X\in\ts{ind}\,\mathcal C$ such that $F_{\lambda}X\in\Gamma_A$ is a
  $G$-stable faithful and generalised standard component of $\Gamma(\ts{mod}\,\mathcal C)$
  with no non-trivial path of the form $X\rightsquigarrow
  \,^gX$. Finally, the map $X\mapsto F_{\lambda}X$ on the vertices of
  $\Gamma_{\mathcal C}$ extends to a Galois covering of translation quivers
  $\Gamma_{\mathcal C}\to \Gamma_A$ with group $G$.
\end{prop}

In order to prove this result, we proceed along the following steps:
\begin{enumerate}
\item[(a)] Any $X\in \Gamma_A$ satisfies $X\simeq
  F_{\lambda}\widetilde{X}$ for some $\widetilde{X}\in \ts{ind}\,\mathcal
  C$ such that $G_{\widetilde{X}}=1$.
\item[(b)] $\ts{rad}^{\infty}(X,Y)=0$ for every
  $X,Y\in \Gamma_{\mathcal C}$.
\item[(c)] $\Gamma_{\mathcal C}$ is a component of $\Gamma(\ts{mod}\, \mathcal C)$.
\end{enumerate}
We prove each step in a separate lemma.
\begin{lem}
  \label{lem2.2}
Let $A$ be connected, weakly shod and not quasi-tilted of canonical
type, $\Gamma_A$ a connecting component of $A$ and
$F\colon \mathcal C \to A$ a Galois covering with group $G$ where $\mathcal C$
is locally bounded. Then for every $X\in\Gamma_A$ there exists
$\widetilde{X}\in \ts{ind}\,\mathcal C$ such that $F_{\lambda}\widetilde{X}\simeq
X$ and $G_{\widetilde{X}}=1$.
\end{lem}
\noindent{\textbf{Proof:}} Note that if $Y=\tau_A^mX$ for some
$m\in\mathbb{Z}$, then the conclusion
holds true for $X$ if and only if it hods true for $Y$. We prove the
lemma by induction on $\ts{rk}(K_0(A))$ and begin with the case where
$A$ is tilted. If $A$ is tilted then $\Gamma_A$ has a complete slice
$\{T_1,\ldots,T_n\}$. By \cite[Cor. 4.5, Prop. 4.6]{lemeur6} and the
above remark, the lemma holds true for $A$. Now assume that $A$ is
not tilted and that the lemma holds true for algebras whose rank of
the Grothendieck group is smaller than $\ts{rk}(K_0(A))$. So $\mathcal
P_A^f\neq\emptyset$. Let
$P_m\in\mathcal P_A^f$ be maximal and $A=B[M]$ the associated
one-point extension. Recall (\cite[Prop. 6.1.40,
Prop. 6.1.41]{lemeur_thesis}) that for any component $B'$ of $B$ the
Galois covering $F\colon \mathcal C\to A$ restricts to a Galois covering
$F^{-1}(B')\to B'$ with group $G$. The conclusion of the lemma
clearly holds true for 
$X=P_m$. Let $B'$ be a component of $B$ and $X$
lie in a connecting component of $B'$. By the induction hypothesis, we
have $X\simeq F'_{\lambda}\widetilde X$ where $\widetilde X\in\ts{ind}\, F^{-1}(B')$ is such that
$G_{\widetilde X}=1$ and $F'\colon F^{-1}(B')\to B'$ is the restriction of
$F$. In particular $X\simeq F_{\lambda}\widetilde X$.
By the above remark and \ref{rem1.6}, the proposition therefore holds
true for 
$A$.\sq

\begin{lem}
\label{lem2.3}
  Keep the notations and hypotheses \ref{lem2.2}. Let
  $\Gamma_{\mathcal C}$ be the full subquiver of $\Gamma(\ts{mod}\, \mathcal C)$ generated by
  the modules $X\in \ts{ind}\,\mathcal C$ such that
  $F_{\lambda}X\in\Gamma_A$. Then:
\begin{enumerate}
\item[(a)]  $\Gamma_{\mathcal C}$ is a (disjoint)
  union of components of $\Gamma(\ts{mod}\, \mathcal C)$.
\item[(b)]$\Gamma_{\mathcal C}$ is faithful, has
  no non trivial path of the form $X\rightsquigarrow\,^gX$ and $\ts{rad}^{\infty}(X,Y)=0$ for every
  $X,Y\in\Gamma_{\mathcal C}$.
\end{enumerate}
\end{lem}
\noindent{\textbf{Proof:}} This follows from \ref{prop1.2}
and the fact that $\Gamma_A$ is faithful, generalised standard,
and has no oriented cycle.\sq

\begin{lem}
  \label{lem2.4}
  Keep the notations and hypotheses of \ref{prop2.1}.  Then $\Gamma_{\mathcal
    C}$ is a component of $\Gamma(\ts{mod}\, \mathcal C)$.
\end{lem}
\noindent{\textbf{Proof:}} Following \cite{happel_reiten_smalo}, we
define the left part $\mathcal L_A$ of $\ts{mod}\,A$ as the full subcategory of
$\ts{ind}\,A$ generated by:
\begin{equation}
 \{M\in \ts{ind}\,A\ |\ \ts{pd}_A\,L\leqslant 1\ \text{for every
    predecessor $L$ of $M$ in $\ts{ind}\,A$}\}\notag
\end{equation}
where $\ts{pd}_A$ is the projective dimension. Let $T$ be the direct sum of
the indecomposable $A$-modules which are either $\ts{Ext}$-injective in $\mathcal
L_A$ or not in $\mathcal L_A$ and projective. Then $T$ is a basic tilting
$A$-module (\cite[4.2,4.4]{assem_coelho_trepode}) and for every $X\in
\Gamma_A$ there exists $m\in\mathbb{Z}$ such that $\tau_A^mX$ is a
direct summand of $T$. Fix an indecomposable decomposition
$T=T_1\oplus\ldots\oplus T_n$ in $\ts{mod}\,A$. By \ref{lem2.2}, there exist
$\widetilde{T}_1,\ldots,\widetilde{T}_n\in
\Gamma_{\mathcal C}$ such that $F_{\lambda}\widetilde{T}_i\simeq T_i$ and
$G_{\widetilde{T}_i}=1$, for
every $i$. Let $\mathcal E$ be the full subcategory of $\ts{ind}\,\mathcal C$
generated by
$\{\,^g\widetilde{T}_i\ |\ i\in\{1,\ldots,n\}\ \text{and}\ g\in
G\}$. So $\mathcal C$ and $\mathcal E$ have equivalent derived categories
(see the proof of \cite[Lem. 4.8]{lemeur6}). In particular $\mathcal E$ is connected. So,
by
\ref{lem2.3}, (b), there is a component $\Gamma$ of
$\Gamma_{\mathcal C}$ which contains $\{\,^g\widetilde{T}_i\}_{i,g}$. We claim
that $\Gamma=\Gamma_{\mathcal C}$. If $X\in \Gamma_{\mathcal
  C}$, then $F_{\lambda}X\in\Gamma_A$ so that
$\tau_A^mF_{\lambda}X\simeq T_i$ for some $i\in\{1,\ldots,n\}$ and
$m\in\mathbb{Z}$. Consequently $\tau_{\mathcal C}^mX\simeq
\,^g\widetilde{T}_i$ for some $g$ and therefore $X\in \Gamma$. Thus
$\Gamma_{\mathcal C}=\Gamma$ is connected.\sq

Now we prove \ref{prop2.1}.

\noindent{\textbf{Proof of \ref{prop2.1}:}} The proposition
is a direct consequence of \ref{lem2.2}, \ref{lem2.3} and \ref{lem2.4}.\sq

\begin{rem}
  \label{rem2.6}
Assume that $\mathcal P_A^f=\emptyset$ and $A$ admits two connecting
components: Its unique postprojective component and its unique
preinjective one. With the hypotheses and notations of \ref{prop2.1}, assume that $\Gamma_A$ is the
postprojective component (or the preinjective component) of $A$. Then it is not difficult to check
that $\Gamma_{\mathcal C}$ is the unique postprojective component (or
the unique preinjective component, respectively) of
$\Gamma(\ts{mod}\,\mathcal C)$.
\end{rem}

\subparagraph{Semi-regular components of the first kind}$\ $

Now we treat the case of semi-regular components containing a
projective or an injective. Most of the work in this
paragraph is based on the following lemma which does not assume $A$
to be weakly shod.
\begin{lem}
  \label{lem2.7}
Let $F\colon \mathcal C\to A$ be a Galois covering with group $G$ where $\mathcal
C$ is locally bounded. Let $\Gamma$ be a component of $\Gamma(\ts{mod}\, A)$ such
that:
\begin{enumerate}
\item[(a)] $\Gamma$ has no multiple arrows and every vertex in $\Gamma$
  is the source of at most two arrows and the target of at most two
  arrows.
\item[(b)] There exists $M_0\in \Gamma$ which is either the source of exactly
  one arrow or the target of exactly one arrow, and which is
  isomorphic to $F_{\lambda}\widetilde{M}_0$ where $\widetilde M_0\in \ts{ind}\,\mathcal C$ is
  such that $G_{\widetilde M_0}=1$.
\end{enumerate}
Then every $X\in\Gamma$ is isomorphic to $F_{\lambda}\widetilde X$ for
some $\widetilde X\in \ts{ind}\,\mathcal C$ such that $G_{\widetilde X}=1$.
\end{lem}
\noindent{\textbf{Proof:}} Let $\mathcal X$ be the set of those modules $X\in\Gamma$ for
which the conclusion of the lemma holds. Therefore $\mathcal X$ contains
$M_0$ and $\mathcal X$ is stable under $\tau_A$ and
$\tau_A^{-1}$ because of \ref{lem1.1}, (d). Assume by absurd that $\mathcal
X\subsetneq\Gamma$. Then by considering an unoriented path in
$\Gamma$ starting from a module $X\in\Gamma\backslash\mathcal X$, ending at
$M_0$ and of minimal length, we have the following (or its dual
treated dually):  There exists an irreducible morphism $u\colon Y\to
X$ with $X\in\mathcal X$, $Y\in\Gamma\backslash\mathcal X$ and such that if $E\to X$ is right minimal almost
split, then either $E=Y$, or $E=Y\oplus Y'$ for some $Y'\in \mathcal X$. We
 prove that $Y\simeq F_{\lambda}\widetilde Y$ for some $\widetilde Y\in \ts{ind}\,\mathcal
C$. For this purpose, we distinguish two cases according to whether $E$
is indecomposable or not. We fix $\widetilde X\in \ts{ind}\,\mathcal C$ such that
$F_{\lambda}\widetilde X\simeq X$ and $G_{\widetilde X}=1$. Assume first that $E=Y$ is
indecomposable. Let $\widetilde u\colon \widetilde Y\to \widetilde X$ be a right minimal
almost split morphism in $\ts{mod}\, \mathcal C$. Thus \ref{lem1.1}, (c), implies
that so is $F_{\lambda}(\widetilde u)\colon F_{\lambda}\widetilde Y\to F_{\lambda}\widetilde
X$. Therefore $F_{\lambda}\widetilde Y\simeq Y$. Now assume that $E=Y\oplus
Y'$ with $Y'\in \mathcal X$. In particular, $Y'\simeq F_{\lambda}\ Y'$ for
some $\widetilde Y'\in \ts{ind}\,\mathcal C$. We thus have a right minimal almost split
morphism $[u,u']\colon Y\oplus Y'\to X$ in $\ts{mod}\, A$. Let $f\colon \widetilde E\to \widetilde X$ be a
right minimal almost split morphism in $\ts{mod}\, \mathcal C$. As above, we deduce that so is
$F_{\lambda}(f)\colon F_{\lambda}\widetilde E\to F_{\lambda}\widetilde X$ in $\ts{mod}\,A$. Therefore
$F_{\lambda}\widetilde E\simeq Y\oplus F_{\lambda}\widetilde Y'$. Applying $F_.$
yields $\bigoplus\limits_{g\in G}\,^g\widetilde E\simeq
F_.Y\oplus\bigoplus\limits_{g\in G}\,^g\widetilde Y'$. Since $\widetilde Y\in \ts{ind}\,\mathcal
C$, we deduce that $^g\widetilde E=\widetilde Y'\oplus\widetilde Y$ for some $g\in G$ and some
$\widetilde Y\in \ts{mod}\,\mathcal C$. Consequently $F_{\lambda}\widetilde E\simeq
F_{\lambda}\widetilde Y'\oplus F_{\lambda}\widetilde Y$ and finally $Y\simeq F_{\lambda}\widetilde
Y$. Hence, in any case, we
have $Y\simeq F_{\lambda}\widetilde Y$ and an irreducible morphism $\widetilde
Y\to \widetilde X$ for some $\widetilde Y\in \ts{ind}\,\mathcal C$. Since $Y\not\in\mathcal X$,
there exists $g\in G\backslash\{1\}$ such that $^g\widetilde
Y\simeq\widetilde Y$. Therefore the morphism $\widetilde Y\to\widetilde X$ defines two irreducible
morphisms $\alpha\colon \widetilde Y\to \widetilde X$ and $\beta\colon \widetilde Y\to \,^g\widetilde
X$. Since $G_{\widetilde X}=1$, and by \ref{lem1.1}, (c), both
$F_{\lambda}(\alpha)\colon F_{\lambda}\widetilde Y\to F_{\lambda}\widetilde X$ and
$F_{\lambda}(\beta)\colon F_{\lambda}\widetilde Y\to F_{\lambda}\,^g\widetilde
X=F_{\lambda}\widetilde X$ are irreducible. On the other hand, $\Gamma$ has no
multiple arrows so there exists an isomorphism $\varphi\colon
F_{\lambda}\widetilde X\xrightarrow{\sim} F_{\lambda}\widetilde X$ such that
$F_{\lambda}(\beta)=\varphi\circ F_{\lambda}(\alpha)$. By the
covering property of $F_{\lambda}$ we have
$\varphi=\sum\limits_{h\in G}F_{\lambda}(\varphi_h)$ with $(\varphi_h)_h\in
\bigoplus\limits_{h\in G} \ts{Hom}_{\mathcal C}(\widetilde Y,\,^h\widetilde X)$. So
$F_{\lambda}(\beta)=\sum\limits_{g\in G}F_{\lambda}(\varphi_h\circ \alpha)$ and
therefore $\beta=\varphi_g\circ \alpha$ because of the covering property of
$F_{\lambda}$. This implies
that $\varphi_g\colon\widetilde X\to \,^g\widetilde X$ is a retraction and therefore an
isomorphism. We get a contradiction because $G_{\widetilde X}=1$.\sq

We apply this lemma to our situation where $A$ is weakly shod and not
quasi-tilted of canonical type.
\begin{prop}
  \label{prop2.8}
Let $A$ be connected, weakly shod and not quasi-tilted of canonical
type, $F\colon \mathcal C\to A$ a Galois covering with group $G$
where $\mathcal C$ is locally bounded and $\Gamma$  a
semi-regular component of $\Gamma(\ts{mod}\, A)$ containing a projective or
an injective. Then for every $X\in\Gamma$ there exists $\widetilde
X\in \ts{ind}\, \mathcal C$ such that $F_{\lambda}\widetilde X\simeq X$ and
$G_{\widetilde X}=1$.
\end{prop}
\noindent{\textbf{Proof:}} It follows from
\cite[6.2]{coelho_lanzilotta} that at least one of the following cases is
satisfied:
\begin{enumerate}
\item[(a)] $\Gamma$ is a postprojective or a preinjective component.
\item[(b)] $\Gamma$ is obtained from a tube or from $\mathbb{Z}A_{\infty}$
  by ray or coray insertions.
\end{enumerate}
In case (a), the proposition follows from: \ref{lem1.1}, (d); the
fact that the $G$-action on $\ts{mod}\,\mathcal C$ commutes with $\tau_{\mathcal
  C}$; and, the fact that
the conclusion of the proposition holds true for indecomposable
projective or injective modules. In case (b), there exists a projective $M_0\in\Gamma$ such
that $\Gamma$ and $M_0$ statisfy the conditions of
\ref{lem2.7}. Whence the proposition.\sq

\begin{rem}
  \label{rem2.9}
Keep the notations and hypotheses of the \ref{prop2.8}.
Let $\widetilde{\Gamma}$ be the full subquiver of $\Gamma(\ts{mod}\, \mathcal C)$ generated
by the
vertices $X\in \ts{ind}\,\mathcal C$ such that $F_{\lambda}X\in
\Gamma$. Then $\widetilde{\Gamma}$ is a union of semi-regular components and contains a
projective or an injective.
\end{rem}

The following example shows that \ref{prop2.8} does not necessarily
hold for regular components, even for tilted algebras.
\begin{ex}
  Let $A$ be the path algebra of the Kronecker quiver
  $\xymatrix{1\ar@<1pt>[r]^a\ar@<-1pt>[r]_b&2}$. It admits a Galois
  covering $F\colon A'\to A$ with group
  $\mathbb{Z}/2\mathbb{Z}=<\sigma>$ where $A'$ is the path algebra of
  the following quiver of type $\widetilde{\mathbb{A}_3}$:
    
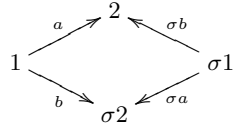
\begin{equation}
      \xymatrix@R=2ex{
&2&\\
1 \ar[ru]^a \ar[rd]_b & &\sigma 1\ar[lu]_{\sigma b} \ar[ld]^{\sigma a}\\
&\sigma 2&
}\notag
    \end{equation}
with $F(x)=F(\sigma x)=x$ for every $x\in\{1,2,a,b\}$. Then the
indecomposable $A$-module
$\xymatrix{k\ar@<1pt>[r]^{\ts{Id}}\ar@<-1pt>[r]_{\ts{id}} & k}$ lying on a homogeneous tube is not of the
first kind with respect to $F$ and, in general, with respect to any
non-trivial connected Galois covering of $A$.
\end{ex}

\section{The Galois covering of $A$ associated to the universal cover
  of the connecting component}$\ $

Let $A$ be weakly shod and not quasi-tilted of canonical type and
$\Gamma_A$ a connecting component. Recall that given a connected Galois covering $F\colon A'\to A$ with
group $G$ there is a component $\Gamma_{A'}$ of
$\Gamma(\ts{mod}\,A')$ and a Galois covering of graphs $\mathcal O(\Gamma_{A'})\to\mathcal
O(\Gamma_A)$ with group $G$ (see \ref{prop2.1} and
\ref{rem1.3}). This Galois covering of graphs is called
\emph{associated to $F$}.
In this section, we prove the following result which is a counter-part
of the work made in \ref{prop2.1}.
\begin{prop}
\label{prop3.1}
  Let $A$ be connected, weakly shod and not quasi-tilted of canonical type,
and  $\Gamma_A$ a connecting component. Then there exists a
  connected Galois covering $ F\colon \widetilde A\to A$ with group 
  the fundamental group $\pi_1(\Gamma_A)$ such that the associated
  Galois covering of graphs $\mathcal O(\Gamma_{\widetilde A})\to \mathcal O(\Gamma_A)$ is the
  universal cover.
\end{prop}

\begin{rem}
\label{rem3.2}
Recall that if $A$ has more than one connecting component, then it has
two of them: The unique preinjective component and the unique
postprojective component. In particular the isomorphism class of
$\pi_1(\Gamma_A)$ does not depend on the
connecting component. 
\end{rem}

Until the end of the section we adopt the hypotheses and notations of
the above proposition. Here is the strategy of its proof. We use an induction on
$\ts{rk}(K_0(A))$.  If $A$ is tilted of type $Q$, then
$\mathcal{O}(\Gamma_A)$ is the underlying graph of
$Q$. So \ref{prop3.1} follows from
\cite[Thm. 1]{lemeur6} in that case. If $A$ is
not tilted, there exists $P_m\in\mathcal P_A^f$ maximal and defining the
one-point extension $A=B[M]$. Then we use \ref{rem1.6} and the Galois covering of $B$
given by the inductive step to construct the desired Galois covering
of $A$.\\

From now on we assume that $A$ is not tilted,
$P_m\in\mathcal P_A^f$ is maximal and $A=B[M]$ is the associated
one-point extension. The extending object is denoted by $x_0\in
A_o$. Also we assume that \ref{prop3.1} holds true for
the components $B_1,\ldots,B_t$  of $B$ ($B=B_1\times\ldots\times
B_t$). Thus for every $i\in\{1,\ldots,t\}$ there is a connected Galois covering $F^{(i)}\colon
\widetilde{B}_i\to B_i$ with group $\pi_1(\Gamma_i)$ equal to the fundamental group
of  the (unique) connecting component $\Gamma_i$ of $B_i$
containing a direct summand of $M$. We write $\widetilde{\Gamma}_i\to\Gamma_i$ for the
universal cover of translation quivers. 
The construction of a connected Galois covering $F\colon \mathcal C\to A$ with group
$\pi_1(\Gamma_A)$ is decomposed into the following steps.
\begin{enumerate}
\item[(a)] A reminder on the universal cover of $\mathcal{O}(\Gamma_A)$.
\item[(b)] The construction of a  Galois covering $F\colon\widetilde B\to B$
   with group $\pi_1(\Gamma_A)$ using  $F^{(1)},\ldots,F^{(t)}$.
\item[(c)] The construction of the locally bounded $k$-category $\widetilde A$ and
  the Galois covering $F\colon \widetilde A\to A$.
\item[(d)] The proof that $\widetilde A$ is connected.
\end{enumerate}

\subparagraph{Reminder: the universal cover 
  $\mathcal{O}(\Gamma_A)$}$\ $

For simplicity, we still denote by
$x_0$ the vertex $(P_m)^{\tau_A}$ of $\mathcal O(\Gamma_A)$ and use it as the
base-point for the computation of the universal cover of
$\mathcal{O}(\Gamma_A)$. Recall
 that the universal cover $p\colon \widetilde{\mathcal O}\to\mathcal O(\Gamma_A)$ is such
 that:
 \begin{enumerate}
 \item[(a)] $\widetilde{\mathcal O}$ is the graph with vertices the homotopy classes $[\Gamma]$ of paths
   $\Gamma\colon x_0\rightsquigarrow x$ in
   $\mathcal{O}(\Gamma_A)$ (where $x$ is any vertex) and such that for
   every edge  $\alpha\colon x-y$ in
  $\mathcal{O}(\Gamma_A)$ and every vertex $[\Gamma]$ in $\widetilde{\mathcal O}$ with
  end-point $x$, there
  is an edge  $\alpha\colon[\Gamma]-[\alpha\Gamma]$ in $\widetilde{\mathcal O}$.
\item[(b)] With the notations of (a), $p$ maps the vertex $[\Gamma]$
  to $x$ and the edge $\alpha\colon[\Gamma]-[\alpha\Gamma]$ to $\alpha\colon x-y$.
 \end{enumerate}

\subparagraph{The Galois covering of  $F\colon\widetilde B\to B$ with group $\pi_1(\Gamma_A)$}$\ $

We construct a Galois covering $F\colon \widetilde B\to B$ with group
$\pi_1(\Gamma_A)$ using $F^{(1)},\ldots,F^{(t)}$. We define $\widetilde B$ as a
disjoint union $\coprod\limits_{i=1}^t\coprod\limits_?\widetilde B_i$ of (infinitely many) copies of $\widetilde B_i$
($i\in\{1,\ldots,t\}$). More precisely, let
$i\in\{1,\ldots,t\}$. Every component $\mathcal U$ of $p^{-1}(\mathcal O(\Gamma_i))$
is simply connected so the restriction $\mathcal U\to \mathcal O(\Gamma_i)$ of $p$
fits into a commutative diagram of graphs:
\begin{equation}
  \xymatrix@=3ex{
\mathcal U \ar[rr]^{\sim} \ar[rd]& & \mathcal O(\widetilde{\Gamma}_i)\ar[ld]\\
&\mathcal O(\Gamma_i)
}\tag{$D_{\mathcal U}$}
\end{equation}
where the horizontal arrow is an isomorphism and the oblique arrow on
the right is induced by $\widetilde{\Gamma}_i\to\Gamma_i$. We then attach to $\widetilde B$ one copy of $\widetilde B_i$ for
each component $\mathcal U$ of $p^{-1}(\mathcal O(\Gamma_i))$. The Galois coverings
$F^{(1)},\ldots,F^{(t)}$ then clearly define a functor $F\colon \widetilde
B\to B$ such that $F$ and $F^{(i)}$ coincide on each copy of $\widetilde B_i$.

Now we endow $\widetilde B$ with a $\pi_1(\Gamma_A)$-action such that $F\circ g=g$
for every $g\in\pi_1(\Gamma_A)$. Let $g\in\pi_1(\Gamma_A)$ and $\widehat{B}_i$ be a
copy of $\widetilde B_i$ in $\widetilde B$. We define the action of $g$ on
$\widehat{B}_i$. Let $\mathcal U$ be the component of $p^{-1}(\mathcal O(\Gamma_i))$
associated to $\widehat{B}_i$. Then $g(\mathcal U)$ is also a component of
$p^{-1}(\mathcal O(\Gamma_i))$ to
which corresponds a copy $\overline{B}_i$ of $\widetilde B_i$ in $\widetilde
B$. Moreover, the graph morphism $g\colon \mathcal U\to g(\mathcal U)$ and the
diagrams ($D_{\mathcal U}$) and ($D_{g(\mathcal U)}$) determine an automorphism
$\mathcal O(\widetilde{\Gamma}_i)\xrightarrow{\sim}\mathcal
O(\widetilde{\Gamma}_i)$ making the following diagram
commute:
\begin{equation}
  \xymatrix@=3ex{
\mathcal O(\widetilde{\Gamma}_i) \ar[rr]^{\sim} \ar[rd] & & \mathcal O(\widetilde{\Gamma}_i) \ar[ld]\\
&\mathcal O(\Gamma_i)&&.
}\notag
\end{equation}
Therefore, the automorphism $\mathcal
O(\widetilde{\Gamma}_i)\xrightarrow{\sim}\mathcal
O(\widetilde{\Gamma}_i)$
extends the map $(X)^{\tau_{\widetilde B_i}}\mapsto
(\,^{\bar{g}}X)^{\tau_{\widetilde
    B_i}}$ associated to some $\bar{g}\in \pi_1(\Gamma_i)$ (see
\ref{rem1.3}). The action of $g$ on $\widehat{B}_i$ is
therefore defined as follows: $g$ maps the component $\widehat{B}_i$
of $\widetilde B$ to the component $\overline{B}_i$ and, as a functor, it acts
like $\bar{g}\colon \widetilde B_i=\widehat{B}_i\xrightarrow{\sim} \widetilde B_i=\overline{B}_i$. This way, we get
a $\pi_1(\Gamma_A)$-action on $\widetilde B$ such that $F\circ g=F$ for every $g\in G$.
\begin{lem}
\label{lem3.4}  
  The $\pi_1(\Gamma_A)$-action on $\widetilde B$ is free, $\widetilde B$ is locally
  bounded and $F\colon\widetilde B\to B$ is a
  Galois covering with group $\pi_1(\Gamma_A)$.
\end{lem}
\noindent{\textbf{Proof:}} Let $x\in\widetilde B_o$ and $g\in\pi_1(\Gamma_A)$ be such
that $gx=x$. We write $\widehat{B}_i$ for the copy of $\widetilde B_i$ in $\widetilde
B$ containing $x$ and $\mathcal U$ for the corresponding component of
$p^{-1}(\mathcal O(\Gamma_i))$. In particular, $g(\mathcal U)=\mathcal U$ and there exists
$g'\in \pi_1(\Gamma_i)$ such that the action of $g$ on $\widehat{B}_i$ is given by
$g'\colon \widetilde B_i=\widehat{B}_i\xrightarrow{\sim} \widehat{B}_i=\widetilde
B_i$. Since $gx=x$, this means that $g'x=x$. So $g'=\ts{Id}_{\widetilde B_i}$
and $g$ is the identity map on $\mathcal U$. Thus, $g$ is the identity on
the universal cover $\mathcal O(\widetilde{\Gamma}_i)$ and therefore on $\widetilde B$. This proves
that the $\pi_1(\Gamma_A)$-action on $\widetilde B$ is free. 

By construction, $\widetilde B$ is locally bounded.

Now we prove that $\pi_1(\Gamma_A)$ acts transitively on $F^{-1}(x)$ for
every $x\in\widetilde B_o$. Let $x,y\in\widetilde B_o$ be such that $Fx=Fy$. By
construction of $F$, there exists $i$ such that $x$ and $y$ lie on
copies $\widehat{B}_i$ and $\overline{B}_i$ of $\widetilde B_i$ in $\widetilde B$, respectively. We
write $\mathcal U$ and $\mathcal V$ for the components of $p^{-1}(\mathcal O(\Gamma_i))$
corresponding to $\widehat{B}_i$ and $\overline{B}_i$,
respectively. So there exists $g\in\pi_1(\Gamma_A)$ such that $g(\mathcal U)=\mathcal
V$. Therefore $gx$ lies on $\overline{B}_i$ and $F(gx)=Fy$. So we may assume that
$\overline{B}_i=\widehat{B}_i$. Using ($\mathcal D_{\mathcal U}$), we
identify the map $\mathcal U\to \mathcal O(\Gamma_i)$ induced by $p$ with the universal
cover $\mathcal O(\widetilde{\Gamma}_i)\to\mathcal O(\Gamma_i)$. Since $F$ coincides with
$F^{(i)}\colon\widetilde B_i\to B_i$ on $\widehat{B}_i$, there exists $g'\in
\pi_1(\Gamma_i)$ such that $g'(x)=y$. Moreover, there exists $g''\in \pi_1(\Gamma_A)$
such that $g''$ and $g'$ coincide on some vertex of $\mathcal U$ (because
$p\colon \widetilde{\mathcal O}\to \mathcal O(\Gamma_A)$ is a Galois
covering with group $\pi_1(\Gamma_A)$) and 
therefore on $\mathcal U$ (because $\mathcal U\to \mathcal
O(\Gamma_i)$ is a Galois covering). We thus have $g''x=y$ with 
$g''\in\pi_1(\Gamma_A)$. This shows the transitivity of
$\pi_1(\Gamma_A)$ on the fibres of 
$F\colon\widetilde B_o\to B_o$.

Therefore $F\colon\widetilde
B\to B$ is, by construction, a covering functor, $\pi_1(\Gamma_A)$ is a group acting freely
on $\widetilde B$ such that $F\circ g=g$ for every $g\in \pi_1(\Gamma_A)$
and $\pi_1(\Gamma_A)$ acts transitively on the fibres of $F\colon \widetilde B_o\to
B_o$. So $F$ is a Galois covering with group $\pi_1(\Gamma_A)$.\sq

\subparagraph{The Galois covering $F\colon\widetilde A \to A$
 with group $\pi_1(\Gamma_A)$}$\ $

Now we extend $F\colon\widetilde B\to B$ to a Galois covering $F\colon\widetilde A\to A$
with group $\pi_1(\Gamma_A)$. Recall that $A=B[M]$. Accordingly
let $\widetilde A$ be the category:
\begin{equation}
  \widetilde A=
  \begin{bmatrix}
    S&\widetilde M\\0&\widetilde B
  \end{bmatrix}\tag{$\star$}
\end{equation}
where $S$ is the category with objects set $S_o=\pi_1(\Gamma_A)\times
\{x_0\}$ and no non-zero morphism except the scalar multiples of the
identity morphisms and $\widetilde M$ is an $S-\widetilde B$-bimodule defined as follows. Fix
an indecomposable decomposition
$M=\bigoplus\limits_{i=1}^t\bigoplus\limits_{j=1}^{n_i}M_{i,j}$ such
that $M_{i,j}\in\ts{ind}\,B_i$ for every $i,j$. Let $i,j$ be such
indices. Then the homotopy class of the edge $x_0-(M_{i,j})^{\tau_A}$ associated to the
inclusion morphism $M_{i,j}\hookrightarrow P_m$ is a vertex in $\widetilde{\mathcal
  O}$ (see \ref{rem1.6}). Also it lies on some component $\mathcal U$ of $p^{-1}(\mathcal O(\Gamma_i))$
to which corresponds a copy $\widehat{B}_i$ of $\widetilde B_i$ in $\widetilde B$. By
\ref{lem2.2}, there exists $\widetilde M_{i,j}\in\ts{ind}\,\widetilde B_i$
such that $F^{(i)}_{\lambda}\widetilde M_{i,j}=M_{i,j}$. We thus consider $\widetilde M_{i,j}$
as an indecomposable $\widetilde B$-module such that $\widetilde
M_{i,j}\in\ts{ind}\,\widehat{B}_i$. In particular we have $F_{\lambda}\widetilde M_{i,j}=M_{i,j}$.
The $S-\widetilde B$-bimodule $\widetilde M$ is then defined as follows:
\begin{equation}
  \begin{array}{crcl}
    \widetilde M\colon & S\times \widetilde B^{op} & \to & \ts{mod}\, k\\
    & ((g,x_0),x) & \mapsto &
    \bigoplus\limits_{i=1}^t\bigoplus\limits_{j=1}^{n_i}\,^g\widetilde
    M_{i,j}(x)\ .
  \end{array}\notag
\end{equation}
The $k$-category $\widetilde A$ is thus completely defined.
Now we extend the $\pi_1(\Gamma_A)$-action on $\widetilde B$ to an action on $\widetilde A$.
We let $\pi_1(\Gamma_A)$ act on $\pi_1(\Gamma_A)\times\{x_0\}$ in the obvious
way. Let $g\in\pi_1(\Gamma_A)$
and $u\in \widetilde M_{i,j}(h^{-1}x)\subseteq \widetilde M((h,x_0),x)$.
We define $g.u$ to be the morphism $u$ viewed as an element of 
$\widetilde M_{i,j}(h^{-1}x)\subseteq \widetilde M(g.(h,x_0),g.x)$.
\begin{lem}
\label{lem3.7}
  The above construction defines a locally bounded $k$-category $\widetilde A$
  endowed with a free $\pi_1(\Gamma_A)$-action.
\end{lem}
\noindent{\textbf{Proof:}} We clearly have defined a $k$-category and
the $\pi_1(\Gamma_A)$-action is well-defined and free because $\pi_1(\Gamma_A)$ acts
freely on $\pi_1(\Gamma_A)\times\{x_0\}$ and on $\widetilde B$. We prove that $\widetilde A$ is
locally bounded. Recall that $\widetilde B$ is locally bounded. Moreover  for
every $g\in \pi_1(\Gamma_A)$ we have 
$\bigoplus\limits_{x\in\widetilde B_o}\widetilde
A((g,x_0),x)=\bigoplus\limits_{x\in \widetilde B_o,i,j}\widetilde
  M_{i,j}(g^{-1}x)=\bigoplus\limits_{x\in B_o}M(x)$ 
because $F_{\lambda}\widetilde
M_{i,j}=M_{i,j}$ for every $i,j$. Thus $\bigoplus\limits_{x\in
    \widetilde B_o}\widetilde A((g,x_0),x)$ is finite dimensional for every $g\in
  \pi_1(\Gamma_A)$. Finally, for every $x\in \widetilde B_o$ we have
  $\bigoplus\limits_{g\in
    \pi_1(\Gamma_A)}\widetilde
  A((g,x_0),x)=\bigoplus\limits_{g\in\pi_1(\Gamma_A),i,j}\widetilde
  M_{i,j}(g^{-1}x)=M(F(x))$.
So $\bigoplus\limits_{g\in
    \pi_1(\Gamma_A)}\widetilde
  A((g,x_0),x)$ is finite dimensional for every $x\in \widetilde B_o$. This
  proves that $\widetilde A$ is locally bounded.
\sq

We extend the Galois covering $F\colon \widetilde B\to B$ to a functor
$F\colon\widetilde A\to A$ as follows:
\begin{enumerate}
\item[(a)] $F((g,x_0))=x_0$ for every $g\in \pi_1(\Gamma_A)$.
\item[(b)] Let $u\in\widetilde M_{i,j}(g^{-1}x)\subseteq \widetilde M((g,x_0),x)$. Then
  $\widetilde M_{i,j}(g^{-1}x)\subseteq \bigoplus\limits_{h\in\pi_1(\Gamma_A)}\widetilde
  M_{i,j}(h^{-1}x)=M_{i,j}(F(x))\subseteq M(F(x))$ (recall that
  $F_{\lambda}\widetilde M_{i,j}=M_{i,j}$). So we set $F(u)=u\in M(F(x))$.
\end{enumerate}

\begin{lem}
\label{lem3.8}
  The above construction defines a Galois covering $F\colon\widetilde
  A\to A$ with group $\pi_1(\Gamma_A)$.
\end{lem}
\noindent{\textbf{Proof:}} 
$F\colon\widetilde A\to A$ is a $k$-linear
functor such that $F\circ g=g$ for every $g\in\pi_1(\Gamma_A)$. Moreover it
is a covering functor
because so is $F\colon \widetilde B\to B$ and 
$F_{\lambda}\widetilde M_{i,j}=M_{i,j}$ for every $i,j$. Finally, the group
$\pi_1(\Gamma_A)$ acts transitively on $F^{-1}(x)$ for every $x\in\widetilde
A_o$. Indeed, this is the case if $x\in\widetilde B_o$ because $F\colon\widetilde B\to
B$ is a Galois covering with group $\pi_1(\Gamma_A)$ and it is clearly the
case if $x=x_0$. So $F$ is a Galois covering with group
$\pi_1(\Gamma_A)$.\sq

\subparagraph{The category $\widetilde A$ is connected}$\ $

We denote by $\widetilde
P_m$ the indecomposable projective $\widetilde A$-module associated to the
object $(1,x_0)$ of $\widetilde A$. Therefore
$\ts{rad}(P_m)=\bigoplus\limits_{i,j}\widetilde M_{i,j}$. We need the
following lemma whose proof follows from the definitions and where
$\widetilde x_0=(x_0,1)$.
\begin{lem}
  \label{lem3.3}
Let $g\in \pi_1(\Gamma_A)$ and $g\widetilde x_0-x_1$ be an edge in $\widetilde{\mathcal O}$. Then
there exist $i,j$ such that $x_1$ is the homotopy class of the edge
 $\alpha\colon x_0-(M_{i,j})^{\tau_A}$ 
in $\mathcal O(\Gamma_A)$ associated to the inclusion $M_{i,j}\hookrightarrow
P_m$. Let $\mathcal U$ be the component of $p^{-1}(\mathcal O(\Gamma_i))$ containing
$x_1$ and $\widehat{B}_i$ the associated copy of $\widetilde B_i$ in $\widetilde
B$. Then $^g\widetilde M_{i,j}\in\ts{ind}\,\widehat{B}_i$ (and $\widetilde M_{i,j}$ is
a direct summand of $\ts{rad}(\,^g\widetilde P_m)$).
\end{lem}

We use \ref{lem3.3} to prove that $\widetilde A$ is connected.
\begin{lem}
  \label{lem3.5}
$\widetilde A$ is connected
\end{lem}
\noindent{\textbf{Proof:}}
It suffices to prove that two indecomposable projective $\widetilde A$-modules lie on the
same component of $\ts{mod}\,\widetilde A$.
Let $g\in\pi_1(\Gamma_A)$. Since $\widetilde O$
is connected, there is a sequence of
edges in $\widetilde O$:
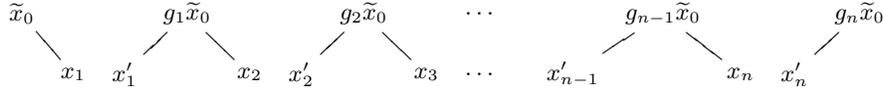
\begin{equation}
  \xymatrix@R=2ex@C=1ex{
\widetilde x_0 \ar@{-}[rd] &&& g_1\widetilde x_0 \ar@{-}[ld] \ar@{-}[rd] &&& g_2\widetilde x_0
\ar@{-}[ld] \ar@{-}[rd]
&& \ldots &
&& g_{n-1}\widetilde x_0 \ar@{-}[ld] \ar@{-}[rd] &&& g_n\widetilde x_0\ar@{-}[ld]\\
&x_1 & x_1' && x_2 & x_2' && x_3 & \ldots &&  x_{n-1}' && x_n& x_n'&
}\notag
\end{equation}
where $g=g_n$ and, for every $j$, the vertices $x_j$ and $x_j'$ lie on
the same component of $p^{-1}(\mathcal O(\Gamma_{i_j}))$ for some $i_j$. By \ref{lem3.3} and because $\widetilde B_1,\ldots,\widetilde B_t$
are connected, the modules $\widetilde P_m$ and $^g\widetilde P_m$ lie on the same
connected component of $\ts{mod}\,\widetilde A$.

Now let $\widetilde P$ be an indecomposable projective
$\widetilde A$-module associated to an object $x\in\widetilde B_o$. So $F_{\lambda}\widetilde P$ is
the indecomposable projective $B$-module associated to $Fx$. Let
$i\in\{1,\ldots,t\}$ be such that $Fx$ is an object of $B_i$. So $x$ is
an object of some copy $\widehat{B}_i$ of $\widetilde B_i$ in $\widetilde B$ and we
let $\mathcal U$ be the associated component of $p^{-1}(\mathcal O(\Gamma_i))$. On the
other hand, we let $\overline{B}_i$ be the copy of $\widetilde B_i$ in $\widetilde B$
such that $\widetilde M_{i,1}\in\ts{ind}\,\overline{B}_i$ and $\mathcal V$ the
associated component of $p^{-1}(\mathcal O(\Gamma_i))$. In particular there
exists $g\in\pi_1(\Gamma_A)$ such that $g(\mathcal V)=\mathcal U$ so that $^g\widetilde
M_{i,1}\in\ts{ind}\,\widehat{B}_i$. Therefore: $\widetilde P$ and $^g\widetilde
M_{i,1}$ lie on the same component of $\ts{mod}\,\widetilde A$ because they
are indecomposable $\widehat{B}_i$-modules and $\widehat{B_i}$ is connected;
$^g\widetilde M_{i,1}$ and $^g\widetilde P_m$ lie on the same component of
$\ts{mod}\,\widetilde A$ because of the inclusion $\widetilde M_{i,j}\hookrightarrow
\widetilde P_m$; and we already proved that so do $\widetilde P_m$ and $^g\widetilde P_m$. This shows
that $\widetilde P$ and $\widetilde P_m$ lie on the same component of
$\ts{mod}\,\widetilde A$. So $\widetilde A$ is connected.\sq

Now we are in position to prove the main result of the section.\\
\noindent\textbf{Proof of \ref{prop3.1}:}
We use an induction on $\ts{rk}(K_0(A))$. If
$A$ is tilted, then the result follows from
\cite[Thm. 1]{lemeur6}. Assume that $A$ is not tilted and that the
conclusion of the proposition holds for algebras $B$ such that
$\ts{rk}(K_0(B))<\ts{rk}(K_0(A))$. Hence there exists a
maximal element $P_m\in\mathcal P_A^f$. Let $A=B[M]$ be the associated
one-point extension. Let $B=B_1\times\ldots\times B_t$ be an
indecomposable decomposition. Then $B_1,\ldots,B_t$ are connected,
weakly shod and not quasi-tilted of canonical type. Let
$\Gamma_1,\ldots,\Gamma_t$ be the connecting components of $B_1,\ldots,B_t$,
respectively, containing a summand of $M$. The induction
hypothesis implies that, for every $i$, there exists a connected Galois covering $F^{(i)}\colon\widetilde
B_i\to B_i$ with group $\pi_1(\Gamma_i)$ whose associated Galois
covering of $\mathcal O(\Gamma_i)$ is the universal cover of graph. By
\ref{lem3.8} and \ref{lem3.5}, there exists a connected Galois
covering $F\colon\widetilde A\to A$ with group $\pi_1(\Gamma_A)$. Let $\mathcal O(\Gamma_{\widetilde
  A})\to\mathcal O(\Gamma_A)$ be the associated Galois covering with group
$\pi_1(\Gamma_A)$. Since $\pi_1(\Gamma_A)$ is free, this
Galois covering is necessarily the universal covering of graphs.\sq

We give some examples to illustrate \ref{prop3.1}. In these
examples we write $P_x$, $I_x$ or $S_x$ for the corresponding indecomposable
projective, indecomposable injective or simple, respectively.
\begin{ex}
\label{ex3.9}
  Let $A$ be the radical square zero algebra with ordinary quiver $Q$
  as follows:
  \begin{equation}
    \xymatrix{
1\ar[r]&2\ar[r]&3\ar[r]&4\ar[r]&5\ar@<2pt>[r]\ar@<-2pt>[r]\ar[r]&6
&.}\notag
  \end{equation}
Let $M=\ts{rad}(P_6)$. Then $A=B[M]$ where $B$ is
the radical square zero algebra with ordinary quiver:
\begin{equation}
  \xymatrix{1\ar[r]&2\ar[r]&3\ar[r]&4\ar[r]&5&.}\notag
\end{equation}
Note that $B$ is of finite representation type and $\Gamma(\ts{mod}\, B)$
is equal to:
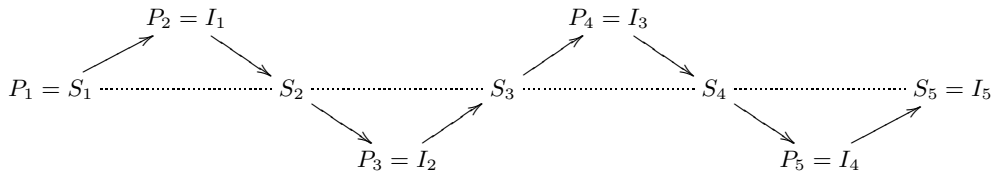
\begin{equation}
  \xymatrix@=1ex{
&& P_2=I_1 \ar[rrdd] && && && && P_4=I_3 \ar[rrdd] && && &&
\\
\\
P_1=S_1 \ar[rruu] \ar@{.}[rrrr] && && S_2 \ar[rrdd] \ar@{.}[rrrr] &&
&& S_3 \ar[rruu] \ar@{.}[rrrr] && && S_4 \ar[rrdd] \ar@{.}[rrrr] && &&
S_5=I_5
\\
\\
&& && && P_3=I_2 \ar[rruu] && && && && P_5=I_4 \ar[rruu] &&&.
}\notag
\end{equation}
The algebra $A$ is wild and weakly shod, it has a unique connecting
component of the following shape:
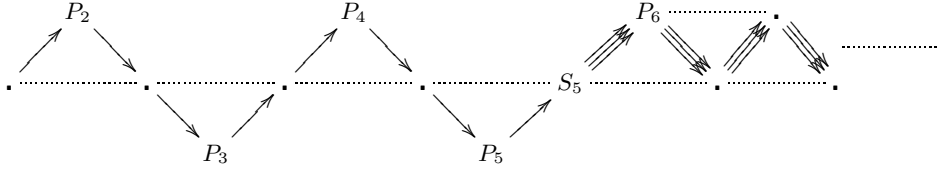
\begin{equation}
  \xymatrix@=1ex{
&& P_2 \ar[rrdd] && && && && P_4 \ar[rrdd] && && && &&
P_6\ar@<2pt>[rrdd] \ar@<-2pt>[rrdd] \ar[rrdd]\ar@{.}[rrrr]&& && \centerdot \ar@<2pt>[rrdd]
\ar@<-2pt>[rrdd]\ar[rrdd]&& &&&&
\\
&&&&&&&&&&&&&&&&&&&&&&&& \ar@{.}[rrrr]&&&&\\
\centerdot \ar[rruu] \ar@{.}[rrrr] && && \centerdot \ar[rrdd] \ar@{.}[rrrr] &&
&& \centerdot \ar[rruu] \ar@{.}[rrrr] && && \centerdot \ar[rrdd] \ar@{.}[rrrr] && &&
S_5 \ar@{.}[rrrr] \ar@<2pt>[rruu] \ar@<-2pt>[rruu]\ar[rruu]&& && \centerdot
\ar@{.}[rrrr] \ar@<2pt>[rruu] \ar@<-2pt>[rruu]\ar[rruu]&& &&\centerdot &&&&
\\
\\
&& && && P_3 \ar[rruu] && && && && P_5 \ar[rruu] &&
 &&
&& && &&&&&&&.
}\notag
\end{equation}
Note that $A$ is not quasi-tilted because the
projective dimension of $S_5$ is equal to $4$. The orbit-graph of the
connecting component of $A$ is equal to:
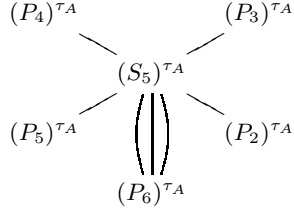
\begin{equation}
  \xymatrix@=2ex{
(P_4)^{\tau_A}&&(P_3)^{\tau_A}
\\
&(S_5)^{\tau_A} \ar@{-}[ru]\ar@{-}[rd]\ar@{-}[lu]\ar@{-}[ld]
\ar@/_/@{-}[dd] \ar@/^/@{-}[dd]\ar@{-}[dd]&\\
(P_5)^{\tau_A}&&(P_2)^{\tau_A}\\
&(P_6)^{\tau_A}&&&.
}\notag
\end{equation}
The fundamental group of this graph is free of rank $2$. So
\ref{prop3.1} implies that $A$ admits a
connected Galois covering with group a free group with rank
$2$. Actually this Galois covering is given by the fundamental group
of the monomial presentation of $A$ (see \cite{martinezvilla_delapena}).
\end{ex}

Recall that weakly shod algebras are particular cases of Laura
algebras. The following example from \cite{coelho_lanzilotta} shows
that \ref{prop3.1} holds for some Laura algebras which are not weakly
shod.
\begin{ex}
 \label{ex3.10}
  (see \cite[2.6]{coelho_lanzilotta}) Let $A$ be the radical square
  zero algebra with ordinary quiver $Q$ as follows:
  
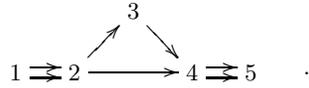
\begin{equation}
\xymatrix@=3ex{
&&3\ar[rd]&&\\
    1 \ar@<2pt>[r]\ar@<-2pt>[r]& 2 \ar[ru]\ar[rr] && 4\ar@<2pt>[r]\ar@<-2pt>[r] &5&.
}\notag
  \end{equation}
Then $A$ is a Laura algebra. The component of $\Gamma(\ts{mod}\, A)$
consist of:
\begin{enumerate}
\item The postprojective components and the homogeneous tubes of the Kronecker
  algebra with quiver $\xymatrix{1 \ar@<2pt>[r]\ar@<-2pt>[r]& 2 }$.
\item The preinjective component and the homogeneous tubes of the Kronecker
  algebra with quiver $\xymatrix{4 \ar@<2pt>[r]\ar@<-2pt>[r]& 5 }$.
\item A unique non semi-regular component of the following shape:
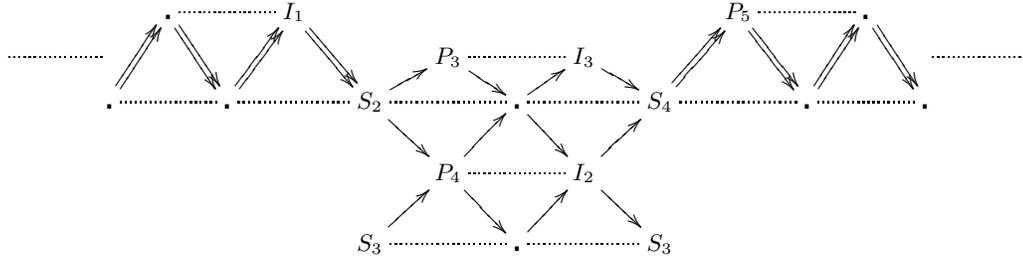
\begin{equation}
  \xymatrix@=1ex{
&&&& &&  \centerdot \ar@<1pt>[rrdd] \ar@<-1pt>[rrdd] \ar@{.}[rrrr]&& && 
I_1 \ar@<1pt>[rrdd] \ar@<-1pt>[rrdd] && && && && && && 
P_5 \ar@<1pt>[rrdd] \ar@<-1pt>[rrdd] \ar@{.}[rrrr] && && 
\centerdot \ar@<1pt>[rrdd] \ar@<-1pt>[rrdd] && &&&&\\
\ar@{.}[rrrr]&& && &&&&&&&&&& P_3 \ar[rrd] \ar@{.}[rrrr] &&&& I_3 \ar[rrd] &&&&&&&&&& \ar@{.}[rrrr] &&&&\\
&&&& \centerdot \ar@<1pt>[rruu] \ar@<-1pt>[rruu] \ar@{.}[rrrr] && && 
 \centerdot \ar@<1pt>[rruu] \ar@<-1pt>[rruu] \ar@{.}[rrrr] && &&
S_2 \ar[rru] \ar[rrdd] \ar@{.}[rrrr]&&  && \centerdot \ar[rru]
\ar[rrdd] \ar@{.}[rrrr]&& 
&& 
S_4 \ar@<1pt>[rruu] \ar@<-1pt>[rruu] \ar@{.}[rrrr] && &&
\centerdot \ar@<1pt>[rruu] \ar@<-1pt>[rruu] \ar@{.}[rrrr] && &&
 \centerdot &&&&\\
\\
&&&&&&&&&&&&&& 
P_4 \ar[rruu] \ar@{.}[rrrr] \ar[rrdd]&& &&
I_2 \ar[rruu] \ar[rrdd] && && && && && &&&&\\
\\
&&&&&& && && && S_3 \ar[rruu] \ar@{.}[rrrr] && && 
\centerdot \ar[rruu] \ar@{.}[rrrr] && && S_3 
}\notag
\end{equation}
where the two copies of the $S_3$ are identified.
\end{enumerate}
In this example, the orbit-graph of the unique non semi-regular
component is the following:
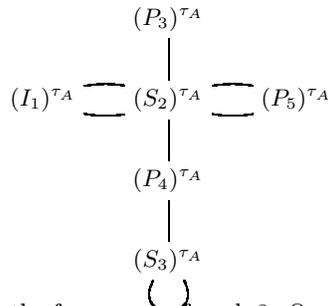
\begin{equation}
  \xymatrix@=4ex{
&(P_3)^{\tau_A} \ar@{-}[d]&\\
(I_1)^{\tau_A}\ar@/^/@{-}[r] \ar@/_/@{-}[r]&(S_2)^{\tau_A} \ar@{-}[d]&(P_5)^{\tau_A}\ar@/^/@{-}[l] \ar@/_/@{-}[l]\\
&(P_4)^{\tau_A}&\\
&(S_3)^{\tau_A} \ar@{-}[u] \ar@{-}@(ld,rd)&
}\notag
\end{equation}
The fundamental group of this graph is the free group of rank
$3$. On the other hand, if one denotes by $(kQ^+)$ for the ideal of
$kQ$ generated by the set of arrows, then the fundamental group of the
natural presentation $kQ/(kQ^+)^2\simeq A$ (in the sense of
\cite{martinezvilla_delapena}) is also isomorphic to the
free group of rank $3$. Hence $A$ admits a connected Galois covering
with group isomorphic to the orbit-graph of the connecting component.
\end{ex}

\section{Proof of Theorem~\ref{thm1} and of Corollary~\ref{cor1}}$\ $

Throughout the section we assume that $A$ is connected and weakly
shod. We prove the first two main results of the text presented in the
introduction.\\
\textbf{Proof of Theorem~\ref{thm1}:} We assume that $A$ is not
quasi-tilted of canonical type. Let $G$ be a group and
$\Gamma_A$ a connecting component of $\Gamma(\ts{mod}\, A)$. If
$F\colon\mathcal C\to A$ is a connected Galois covering then
\ref{prop2.1} yields a Galois
covering of translation quivers with group $G$ of
$\Gamma_A$. Conversely, let $\Gamma'\to
\Gamma_A$ be a Galois covering of translation quivers with group $G$.  Therefore $G\simeq
\pi_1(\Gamma_A)/N$ for some normal subgroup
$N\vartriangleleft \pi_1(\Gamma_A)$ (\cite[1.4]{bongartz_gabriel}). On the other hand, \ref{prop3.1}
yields a connected Galois covering $\widetilde A\to A$ with group
$\pi_1(\Gamma_A)$. Factoring out by $N$ yields a connected Galois covering
$\widetilde A/N\to A$ with group $G$.\sq

Now we turn to the proof of Corollary~\ref{cor1}. We need the three
following lemmas. The first one follows directly from
 Theorem~\ref{thm1} so we omit the proof.
\begin{lem}
\label{lem4.1}
  Assume that $A$ is not quasi-tilted of canonical type. Let $\Gamma_A$ be
  a connecting component of $A$. Then the
  following conditions are equivalent:
  \begin{enumerate}
  \item[(a)] $A$ is simply connected,
  \item[(b)] The orbit-graph $\mathcal{O}(\Gamma_A)$ is a tree.
  \item[(c)]  $\Gamma_A$ is simply connected.
  \end{enumerate}
\end{lem}

The following lemma expresses the simple connectedness of $A=B[M]$ in
terms of the simple connectedness of the components of
$B$. In the case where $A$ is tame weakly shod, the necessity was proved in
\cite[Lem. 5.1]{assem_lanzilotta}. We recall that if $A$ is connected
and $x_0\in A_o$ is
the extension object in $A=B[M]$, then $x_0$ is called \emph{separated}
if $M$ has exactly as many indecomposable summands as the number of
components of $B$ (that is, $M$ restricts to an indecomposable module
on each component of $B$).
\begin{lem}
\label{lem4.2}
  Assume that $\mathcal P_A^f\neq\emptyset$. Let $P_m\in\mathcal P_A^f$ be
  maximal, $A=B[M]$ the associated one-point extension and
  $x_0\in A_o$ the extending object. Then $A$ is simply connected
  if and only if the two following conditions are satisfied:
  \begin{enumerate}
  \item[(a)] $B$ is a product of simply connected algebras,
\item[(b)] $x_0$ is separating (that is, $M$ is multiplicity-free).
  \end{enumerate}
\end{lem}
\noindent{\textbf{Proof:}} By \cite[4.5, 4.8]{coelho_lanzilotta}, $B$
is a product of connected, weakly shod and not quasi-tilted of
canonical type algebras. Assume that $A$ is simply connected. By
\cite[2.6]{assem_delapena}, the object 
$x_0$ is separating. Let $B'$ be
a connected component of $B$. Since $A$ is connected, $M$ admits an
indecomposable summand lying on $\ts{ind}\,B'$. By
\ref{rem1.6} and because the orbit-graph of any connecting component
of $A$ is simply connected, the orbit-graph of any connecting
component of $B'$ is simply connected. So $B'$ is simply connected by
Theorem~\ref{thm1}. Conversely assume that $x_0$ is separating and $B$ is
a product of simply connected algebras. By Theorem~\ref{thm1}, for
every component $B'$ of $B$, the orbit-graph of any connecting
component of $B'$ is a tree. By
\ref{rem1.6} and because $x_0$ is separating, we deduce that
the orbit-graph of any connecting component of $A$ is a tree. By
Theorem~\ref{thm1}, this implies that $A$ is simply connected.\sq

Finally, we recall the following lemma which was proved in
\cite[2.5]{assem_lanzilotta}.
\begin{lem}
\label{lem4.3}
   Under the hypothesis and notations of
  \ref{lem4.2}, the following conditions are equivalent:
  \begin{enumerate}
  \item[(a)] $\ts{HH}^1(A)=0$.
  \item[(b)] $\ts{HH}^1(B)=0$ and $x_0$ is separating.
  \end{enumerate}
\end{lem}

Now we can prove Corollary~\ref{cor1}.

\noindent{\textbf{Proof of Corollary~\ref{cor1}}:} We use an induction
on $\ts{rk}(K_0(A))$. By 
\cite[Thm. 1]{lemeur6}, the corollary holds true if $A$ is
tilted. So we assume that $A$ is not quasi-tilted and
 the corollary holds true for
 algebras $B$ such that
$\ts{rk}(K_0(B))<\ts{rk}(K_0(A))$. Since $\mathcal P_A^f\neq\emptyset$, there exists $P_m\in\mathcal
P_A^f$ maximal. Let $A=B[M]$ be the associated one-point
extension. Using the induction hypothesis applied to the components of
$B$ and using \ref{lem4.1}, \ref{lem4.2} and \ref{lem4.3}, we deduce that $A$ is simply connected if and
only if $\ts{HH}^1(A)=0$. On the other hand, Theorem~\ref{thm1} shows that $A$ is simply
connected if and only if $\mathcal{O}(\Gamma_A)$ is a tree.\sq

We finish this section with an example to illustrate
Corollary~\ref{cor1}
\begin{ex}
\label{ex4.4}
  Let $A$ be as in \ref{ex3.9}. Then $A$ is not simply
  connected and neither is the orbit-graph of its connecting component. On the other hand, a straightforward
  computation shows that
   $\ts{dim}\,\ts{HH}^0(A)=1$,
  $\ts{dim}\, \ts{HH}^1(A)=3$ and
  $\ts{dim}\, \ts{HH}^i(A)=0$ if $i\geqslant 2$.
\end{ex}

\section{The class of weakly shod algebras is stable under finite
  Galois coverings and under quotients}$\ $

In this section we prove Theorem~\ref{thm2}. At first we study the
implications of Theorem~\ref{thm2} in the more general setting of Galois
coverings with non necessarily finite groups.
\begin{lem}
\label{lem5.1}
  Let $F\colon\mathcal C\to A$ be a connected Galois covering with group
  $G$. If $A$ is weakly shod and not quasi-tilted, then
    $\Gamma(\ts{mod}\, \mathcal C)$ has a unique non semi-regular component
    $\Gamma_{\mathcal C}$. Moreover it is faithful, generalised standard and
    has no non
    trivial path
    of the form $X\rightsquigarrow\,^gX$ with $X\in\Gamma_{\mathcal C}$
    and $g\in G$.
\end{lem}
\noindent{\textbf{Proof:}} Let $\Gamma_A$ be the connecting
component of $A$. Let $\Gamma_{\mathcal C}$ be as in
\ref{prop2.1}.  We
only need
to prove that $\Gamma_{\mathcal C}$ is the unique non semi-regular
component of $\Gamma(\ts{mod}\, \mathcal C)$. Note that $\Gamma_{\mathcal C}$ contains
both a projective and an injetive because so does $\Gamma_A$. Let $P\in \ts{ind}\,\mathcal
C\backslash\Gamma_{\mathcal C}$ be projective. Then $F_{\lambda}P\in \ts{ind}\,A\backslash\Gamma_A$
is projective and therefore lies on a semi-regular component of
$\Gamma(\ts{mod}\, A)$. By \ref{rem2.9}, so does $P$. Whence
the lemma.\sq

The preceding lemma has a converse under the additional assumption
that the group $G$ acts freely on the indecomposable modules lying on
$\Gamma_{\mathcal C}$. This last condition is
always verified when $G$ is torsion-free.
\begin{lem}
\label{lem5.2} 
  Let $F\colon\mathcal C\to A$ be a connected Galois covering with group
  $G$. Assume that $\Gamma(\ts{mod}\, \mathcal C)$ has a unique non semi-regular component
    $\Gamma_{\mathcal C}$ and that the following conditions are
    satisfied:
\begin{enumerate}
\item[(a)] $\Gamma_{\mathcal C}$ is faithful and generalised standard.
\item[(b)] $\Gamma_{\mathcal C}$ has no non
    trivial path
    of the form $X\rightsquigarrow\,^gX$.
\item[(c)] $G_X=1$ for every $X\in\Gamma_{\mathcal C}$.
\end{enumerate}
Then $A$ is weakly shod.
\end{lem}
\noindent{\textbf{Proof:}} Note that $\Gamma_{\mathcal C}$ is $G$-stable because
of its uniqueness. If follows from the arguments presented in
the proof of \cite[3.6]{gabriel} that there is a component $\Gamma$ of
$\Gamma(\ts{mod}\, A)$ such that $\Gamma=\{F_{\lambda}X\ |\ X\in\Gamma_{\mathcal
  C}\}$. Also the map $X\mapsto F_{\lambda}X$ extends to a
Galois covering of translation quivers $\Gamma_{\mathcal C}\to \Gamma$
with group $G$. In particular $\Gamma$ is non semi-regular. Moreover
\ref{prop1.2} implies that $\Gamma$ is faithful,
generalised standard and has no oriented cycles. Therefore $A$ is weakly shod.\sq

Now we prove the equivalences of Theorem~\ref{thm2}. Part of the tilted case was
treated in \cite[Rem. 4.10]{lemeur6}. We recall it for convenience.
\begin{prop}
\label{prop6.3}
  Let $F\colon A'\to A$ be a connected Galois covering with finite
  group $G$. Then $A'$ is tilted if $A$ is tilted.
\end{prop}

Now we prove the equivalence of Theorem~\ref{thm2} in the quasi-tilted case.
\begin{prop}
\label{prop6.4}
  Let $F\colon A'\to A$ be a connected Galois covering with finite
  group $G$. Then $A'$ is quasi-tilted if and only if $A$ is quasi-tilted.
\end{prop}
\noindent{\textbf{Proof:}} Recall that $\mathcal L_A$ denotes the
left part of $A$. We use the following description of $\mathcal L_A$
(\cite[Thm. 1.1]{assem_lanzilotta_redondo}):
\begin{equation}
  \mathcal L_A=\{M\in \ts{ind}\, A\ |\
\text{$\ts{pd}_A(L)\leqslant 1$ for every $L\in \ts{ind}\, A$ such that $\ts{Hom}_A(L,M)\neq
  0$}\}\ .\notag
\end{equation}
Also, by (\cite[II Thm. 1.14, II
Thm. 2.3]{happel_reiten_smalo}),  the following conditions are
equivalent for any algebra $A$:
\begin{enumerate}
\item[(a)] $A$ is quasi-tilted.
\item[(b)] $A$ has global dimension at most $2$ and
  $\ts{id}_A(X)\leqslant 1$ or $\ts{pd}_A(X)\leqslant 1$ for every $X\in
  \ts{ind}\, A$.
\item[(c)] $\mathcal L_A$ contains all the indecomposable projective $A$-modules.
\end{enumerate}
Assume that $A$ is quasi-tilted. Let $u\colon X\to P$ be a non-zero morphism of
$A'$-modules
with $X,P\in \ts{ind}\, A'$ and $P$ projective. So $F_{\lambda}(u)\colon
F_{\lambda}X\to F_{\lambda}P$ is non zero and $F_{\lambda}P$ is
indecomposable projective. Fix an indecomposable decomposition $F_{\lambda}X=X_1\oplus\ldots\oplus
X_r$ in $\ts{mod}\,A$. So the restriction $X_i\to F_{\lambda}P$ of
$F_{\lambda}(u)$ is non-zero for some $i$. Since $A$ is quasi-tilted, we
have $F_{\lambda}P\in \mathcal L_A$ and therefore $\ts{pd}_A(X_i)\leqslant
1$. On the other hand, $F_.F_{\lambda}X=\bigoplus\limits_{g\in G}\,^gX$,
$F_.F_{\lambda}X=F_.X_1\oplus\ldots\oplus F_.X_r$ and the projective
dimension is unchanged under $F_.$, $F_{\lambda}$ and under the action
of $G$. Consequently $\ts{pd}_{A'}(X)=\ts{pd}_{A'}(F_.X_i)=\ts{pd}_A(X_i)\leqslant
1$. So $P\in \mathcal L_{A'}$. Thus, $A'$ is quasi-tilted.

Conversely, assume that $A'$ is quasi-tilted. In particular, $A$ and
$A'$ have the same global dimension, that is, at most $2$. Let $X\in
\ts{ind}\, A$. Since $G$ is finite, $F_.X\in \ts{mod}\, A'$. Fix an
indecomposable decomposition
$F_.X=X_1\oplus\ldots \oplus X_r$ in $\ts{mod}\,A'$. We claim that
$X_1,\ldots,X_r$ have the same projective dimension. Indeed, let
$d=\ts{pd}_{A'}(X_1)$ and $I=\{i\in\{1,\ldots,r\}\ |\
  \ts{pd}_{A'}(X_i))=d\}$. Then $F_.X=L\oplus
M$ where $L=\bigoplus\limits_{i\in I} X_i$ and
$M=\bigoplus\limits_{i\in I^c} X_i$. Since the $G$-action on
$\ts{mod}\,A'$ preserves the projective dimension, we have $^gL=L$ and
$^gM=M$ for every $g\in G$. By \cite[1.2]{dowbor_skowronski}, we
deduce that there exist $Y,Z\in\ts{mod}\,A$ such that $X=Y\oplus Z$,
$L=F_.Y$ and $M=F_.Z$. Since $X$ is indecomposable and $I\neq\emptyset$, we have 
$Z=0$ and, therefore, $I=\{1,\ldots,r\}$. Thus
$\ts{pd}_{A'}(X_i)=\ts{pd}_{A'}(X_j)=\ts{pd}_A(X)$ (and, dually,
$\ts{id}_{A'}(X_i)=\ts{id}_{A'}(X_j)=\ts{id}_A(X)$) for every
$i,j$. Since $A'$ is quasi-tilted, we infer that
$\ts{pd}_A(X)\leqslant 1$ of $\ts{id}_A(X)\leqslant 1$. This proves
that $A$ is quasi-tilted.\sq
 
Now we end the proof of Theorem~\ref{thm2}.\\
\noindent{\textbf{Proof of Theorem~\ref{thm2}:}}
The necessity in (a) follows from \ref{prop6.3} and (b) was proved in
\ref{prop6.4}. We prove (c) and may assume that neither $A$ nor $A'$
is 
quasi-tilted.
 Assume that $A$ is weakly shod and not quasi-tilted. Then
\ref{lem5.1} implies that $\Gamma(\ts{mod}\, A')$ has a unique non
semi-regular component which is moreover faithful, generalised
standard and has no oriented cycle. Therefore $A'$ is weakly
shod. This proves the necessity in (c).
 From now on, we assume that $A'$ is weakly shod and not
quasi-tilted of canonical type. We prove that $A$ is weakly shod.
In
view of \ref{lem5.2}, we
need the following result.
\begin{lem}
 \label{lem5.3}  Assume that $A'$ is weakly shod and not quasi-tilted
 of canonical type. We have
 $G_X=1$ for every indecomposable $A'$-module $X$ lying on a
 connecting component of $\Gamma(\ts{mod}\, A')$.
\end{lem}
\noindent{\textbf{Proof of \ref{lem5.3}:}} The conclusion of the lemma holds true
 for any indecomposable projective or injective $A'$-module. So does
 it for non-stable modules because $\tau_{A'}$ commutes with the $G$-action. Let $\Gamma_{A'}$ be a
connecting component of $A'$ and $X\in\Gamma_{A'}$ be
stable. We still write $\mathcal L_{A'}$ for
the left part of $\ts{mod}\, A'$ and we write $\mathcal R_{A'}$ for the right part of
$\ts{mod}\, A'$, defined dually. Since $A'$ is weakly shod, the set
$\ts{ind}\,A'\backslash\left(\mathcal L_{A'}\cup \mathcal R_{A'}\right)$ is finite, 
contained in $\Gamma_{A'}$, and has no periodic module. Therefore there exists $n\in\mathbb{Z}$ such
that $\tau_{A'}^nX\in \Gamma_{A'}\cap \left(\mathcal L_{A'}\cup\mathcal R_{A'}\right)$. Assume for example that
$X'=\tau_{A'}^nX\in\Gamma_{A'}\cap \mathcal L_{A'}$ (the remaining case is dealt with
dually). Let $e$ be the sum of the primitive idempotents $e'$ of $A'$ such
that $e'A'\in \mathcal L_{A'}$ and let $B'=eA'e$. Therefore $B'$ is a full
convex subcategory of $A'$, it is a product of tilted algebras,
$X'\in\ts{ind}\,B'$ (see \cite{assem_coelho_trepode}) and $B'$ is
stable under $G$ because so is
$\mathcal L_{A'}$.
In particular, $F$ restricts to a Galois
covering $F'\colon B'\to B$ with group $G$, where $B:=F(B')$.
In order to prove that $G_X=1$ we prove that $G_{X'}=1$. 
By absurd assume that there is $g\in G\backslash \{1\}$ such that
$^gX'\not\simeq X'$. After replacing $g$ by some adequate power, we
assume that $g$ is of prime order $p$. The quotient $\pi\colon
B\to B/\left<g\right>$ is a Galois covering with group
$\left<g\right>\simeq\mathbb{Z}/p\mathbb{Z}$. Therefore
  $\ts{Ext}^1_{B/\left<g\right>}(\pi_{\lambda}X',\pi_{\lambda}X')\simeq
  \bigoplus\limits_{j=0}^{p-1}\ts{Ext}^1_B(X',\,^{g^j}X')=0$
because of \cite[2.1]{lemeur6}, the isomorphism $^gX'\simeq X'$ and
the equality $\ts{Ext}^1_B(X',X')=0$. In order to get a contradiction we
first prove that $\pi_{\lambda}X'$ is indecomposable. Fix an
indecomposable decomposition $\pi_{\lambda}X'=M_1\oplus\ldots \oplus M_l$ in
$\ts{mod}(B/\left<g\right>)$. Hence
$\ts{Ext}^1_{B/\left<g\right>}(M_i,M_i)=0$ for all $i$. We claim that
$M_i$ lies in the image of $\pi_{\lambda}$ for all $i$. Indeed, we
distinguish two cases according to whether $\ts{car}(k)=p$ or
$\ts{car}(k)\neq p$. If $\ts{car}(k)=p$ then the claim follows from
\cite[Lem. 6.1]{lemeur5}. If $\ts{car}(k)\neq p$, then
$B/\left<g\right>$ is Morita equivalent to the skew-group algebra
$B[\left<g\right>]$ (\cite[Thm. 2.8]{cibils_marcos}) and
$B[\left<g\right>]$ is tilted
(\cite[Thm. 1.2, (g)]{assem_lanzilotta_redondo}). Therefore
$B/\left<g\right>$ is tilted and the claim follows from
\cite[Prop. 4.6]{lemeur6}. Thus, in all cases,
$M_i\simeq\pi_{\lambda}M_i'$ for some
$M_i'\in\ts{mod}\,B$ (necessarily indecomposable). So
$\bigoplus\limits_{j=0}^{p-1} \,^{g^j}M_i'\simeq \pi_.M_i$ is a
summand of
$\pi_.\pi_{\lambda}X'\simeq \bigoplus\limits_{j=0}^{p-1}\,^{g^j}X'\simeq
\bigoplus\limits_{j=0}^{p-1}X'$. We
thus have $M_i'\simeq X'$ for all $i$, whereas
$\pi_{\lambda}X'=\pi_{\lambda}M_1'\oplus\ldots\oplus
\pi_{\lambda}M_l'$. this proves that $\pi_{\lambda}X'$ is
indecomposable. The contradiction is therefore the following. On the
one hand,
$\ts{Ext}_{B/\left<g\right>}^1(\pi_{\lambda}X',\pi_{\lambda}X')=0$,
$\pi_{\lambda}X'\in\ts{ind}(B/\left<g\right>)$ and $B/\left<g\right>$
is a product of quasi-tilted algebras (because $B$ is a product of
tilted algebras and by \ref{prop6.4}), which imply that
$\ts{End}_{B/\left<g\right>}(\pi_{\lambda}X')\simeq k$. On the other
hand,
$\ts{End}_{B/\left<g\right>}(\pi_{\lambda}X') \simeq
\bigoplus\limits_{j=0}^{p-1}\ts{Hom}_B(X',\,^{g^j}X') \simeq
\bigoplus\limits_{j=0}^{p-1}\ts{End}_{B}(X')$ as $k$-vector
spaces. This is absurd. So $G_{X'}=1$ and therefore $G_X=1$.
\sq

Now we can prove that $A$ is weakly shod by applying
\ref{lem5.2}. As remarked in the proof of
\ref{lem5.3}, a non trivial path in $\ts{ind}\,A$ of the form
$X\rightsquigarrow\,^gX$ with $X\in \Gamma_{A'}$ gives rise to a non trivial path
$X\rightsquigarrow X$ in $\ts{ind}\,A'$ which is impossible because
$A'$ is weakly shod. Therefore all the hypotheses of \ref{lem5.2} are
satisfied and
$A$ is weakly shod. This proves (c).

It only remains to prove the necessity in (a). We assume that $A'$ is
tilted and prove that so is $A$. Let $\Gamma_{A'}$ be a 
connecting component of $\Gamma(\ts{mod}\, A')$. It admits a complete
slice $\Sigma'$. Clearly, $\Gamma_{A'}$ is $G$-stable whatever the number
of connecting components of $A'$ is (one or two).
By \ref{lem5.2}, \ref{lem5.3} and \cite[3.6]{gabriel}, there exists a component
$\Gamma$ of $\Gamma(\ts{mod}\, A)$ such that $\Gamma=\{F_{\lambda}X\ |\ X\in\Gamma_{A'}\}$.
Moreover there is a Galois covering of translation quivers $\Gamma_{A'}\to \Gamma$ with
group $G$ extending the map $X\mapsto F_{\lambda}X$.
We prove that $\Gamma$ has a complete slice. For this purpose we use
the following lemma.
\begin{lem}
  \label{lem5.4}
$^gX\in \Sigma'$ for every $g\in G,X\in\Sigma'$.
\end{lem}
\noindent{\textbf{Proof of \ref{lem5.4}:}} Let $g\in G$ and
write $\Sigma'=\{X_1,\ldots,X_n\}$. So there exist a permutation
$i\mapsto g.i$ of $\{1,\ldots,n\}$ and integers $l_1,\ldots,l_n$ such that
 $^gX_i=\tau_{A'}^{l_i}X_{g.i}$ for every $i$. Clearly, the modules
 $^gX_1,\ldots,^gX_n$ form a complete slice $g(\Sigma')$ in $\Gamma_{A'}$. This
 implies that $l_1=l_2=\ldots=l_n$. We write $l=l_1$. Therefore
 $g(\Sigma')=\tau_{A'}^l(\Sigma')$. On the other hand, $g$ has finite
 order and $\Gamma_{A'}$ has no oriented cycles. So $l=0$ and $g(\Sigma')=\Sigma'$.\sq

Let $\Sigma$ be the full subquiver of $\Gamma$ generated by
$\{F_{\lambda}X\ |\ X\in\Sigma'\}$. Hence $\Sigma$ is convex in
$\Gamma$, has no oriented cycle and intersects each $\tau_A$-orbit of $\Gamma$
exactly once because $\Sigma'$ is a $G$-stable complete
slice in $\Gamma_{A'}$.
Moreover, the arguments used in the proof of \ref{prop1.2}
show that $\Sigma$ is faithful because so is $\Sigma'$. Finally,
given $X,Y\in\Sigma'$, we have
$\ts{Hom}_A(F_{\lambda}X,\tau_AF_{\lambda}Y)\simeq \bigoplus\limits_{g\in G}
\ts{Hom}_{A'}(X,\tau_{A'}\,^gY)=0$ because of the covering property of
$F_{\lambda}$, \ref{lem1.1} (d) and the fact that $\Sigma'$ is a $G$-stable
slice in $\Gamma_{A'}$. Thus $\Sigma$ is a complete slice and $A$ is tilted with
$\Gamma$ as a connecting component. This proves the sufficiency (a) and finishes the
proof of Theorem~\ref{thm2}.
\sq

\begin{rem}
The reader may find similar equivalences to those of
Theorem~\ref{thm2} about skew-group algebras (instead of Galois
coverings) under the additional assumption that $\ts{car}(k)$ does not
divide the order of the group $G$ (see \cite{assem_lanzilotta_redondo}).
\end{rem}

\section*{Acknowledgements}

The author gratefully acknowledges Ibrahim Assem for his
encouragements and useful comments.

\bibliographystyle{plain}
\bibliography{biblio}

\noindent Patrick Le Meur\\
\textit{e-mail:} Patrick.LeMeur@cmla.ens-cachan.fr\\
\textit{address:} CMLA, ENS Cachan, CNRS, UniverSud, 61 Avenue du President Wilson, F-94230 Cachan

\end{document}